\documentclass[final,1p,times]{elsarticle}

\usepackage{amssymb}
 \usepackage{amsthm}
\usepackage{amscd}
\usepackage{amsmath}
\usepackage{amsfonts}
\usepackage{amssymb}
\usepackage{graphicx}
\newtheorem{theorem}{Theorem}
\usepackage{mathrsfs}
\usepackage{titletoc}

\newcommand{\la}{\Delta}

\newcommand{\ra}{\rightarrow}
\newcommand{\p}{\partial}
\newcommand{\f}{\frac}

\newcommand{\be}{\begin{equation}}
\renewcommand{\ra}{\rightarrow}
\newcommand{\ee}{\end{equation}}
\newcommand{\bea}{\begin{eqnarray}}
\newcommand{\eea}{\end{eqnarray}}
\newcommand{\bna}{\begin{eqnarray*}}
\newcommand{\ena}{\end{eqnarray*}}

\renewcommand{\le}{\left}
\newcommand{\ri}{\right}

\journal{***}

\begin{document}

\begin{frontmatter}
\title{Extremal functions for Trudinger-Moser inequalities of Adimurthi-Druet type in dimension two}

\author{Yunyan Yang}
 \ead{yunyanyang@ruc.edu.cn}
\address{ Department of Mathematics,
Renmin University of China, Beijing 100872, P. R. China}

\begin{abstract}
  Combining Carleson-Chang's result \cite{CC} with blow-up analysis, we prove existence of extremal functions for certain
  Trudinger-Moser inequalities in dimension two. This kind of inequality was originally proposed by
  Adimurthi and O. Druet \cite{A-D}, extended by the author to high dimensional case and Riemannian surface case
  \cite{Yang-JFA,Yang-Tran}, generalized by C. Tintarev to wider cases including singular form
  \cite{Tint} and by M. de Souza and J. M. do \'O \cite{de-dO} to the whole Euclidean space
 $\mathbb{R}^2$. In addition to the Euclidean case, we also consider the Riemannian surface case. The results in the current
  paper complement that of L. Carleson and A. Chang \cite{CC}, M. Struwe \cite{Struwe},
  M. Flucher \cite{Flucher}, K. Lin \cite{Lin}, and Adimurthi-Druet \cite{A-D}, our previous ones \cite{Yang-Tran,Lu-Yang}, and
  part of C. Tintarev \cite{Tint}.
\end{abstract}


\end{frontmatter}

\titlecontents{section}[0mm]
                       {\vspace{.2\baselineskip}}
                       {\thecontentslabel~\hspace{.5em}}
                        {}
                        {\dotfill\contentspage[{\makebox[0pt][r]{\thecontentspage}}]}
\titlecontents{subsection}[3mm]
                       {\vspace{.2\baselineskip}}
                       {\thecontentslabel~\hspace{.5em}}
                        {}
                       {\dotfill\contentspage[{\makebox[0pt][r]{\thecontentspage}}]}

\setcounter{tocdepth}{2}


\section{Introduction}
Let $\Omega$ be a smooth bounded domain in $\mathbb{R}^2$ and $W_0^{1,2}(\Omega)$ be the usual Sobolev
space. The classical Trudinger-Moser inequality \cite{24,19,17,22,14} says
\be\label{Tr}\sup_{u\in W_0^{1,2}(\Omega),\,\|\nabla u\|_2\leq 1}\int_\Omega e^{4\pi u^2}dx<\infty.\ee
Here and throughout this paper we denote the $L^p$-norm by $\|\cdot\|_p$. This inequality is sharp in the sense that
for any $\alpha>4\pi$, the integrals in (\ref{Tr}) are still finite but the supremum is infinite. Let
$u_k\in W_0^{1,2}(\Omega)$ be such that $\|\nabla u_k\|_2=1$ and $u_k\rightharpoonup u$ weakly in
$W_0^{1,2}(\Omega)$. Then P. L. Lions \cite{Lions} proved that for any $p<1/(1-\|\nabla u\|_2^2)$, there holds
\be\label{Lion}\limsup_{k\ra\infty}\int_\Omega e^{4\pi pu_k^2}dx<\infty.\ee
This inequality gives more information than the Trudinger-Moser inequality (\ref{Tr}) in case $u\not\equiv 0$.
While in case $u\equiv 0$, it is weaker than (\ref{Tr}). However Adimurthi and O. Druet \cite{A-D} proved that for any
$\alpha$, $0\leq\alpha<\lambda_1(\Omega)$,
\be\label{a-d}\sup_{u\in W_0^{1,2}(\Omega),\,\|\nabla u\|_2\leq 1}\int_\Omega e^{4\pi u^2(1+\alpha\|u\|_2^2)}dx<\infty,\ee
and that the supremum is infinity when $\alpha\geq \lambda_1(\Omega)$,
where $\lambda_1(\Omega)$ is the first eigenvalue of the Laplacian operator with respect to Dirichlet boundary condition.
For any sequence of functions $u_k\in W_0^{1,2}(\Omega)$ with $\|\nabla u_k\|_2=1$
and $u_k\rightharpoonup u$ weakly in $W_0^{1,2}(\Omega)$, if $u\not\equiv 0$, it then follows from (\ref{a-d}) that
for any $\alpha$, $0\leq\alpha<\lambda_1(\Omega)$,
\be\label{a-d-k}\limsup_{k\ra\infty}\int_\Omega e^{4\pi u_k^2(1+\alpha \|u_k\|_2^2)}dx<\infty.\ee
Note that $1+\alpha \|u_k\|_2^2<1+\|\nabla u\|_2^2<1/(1-\|\nabla u\|_2^2)$ for sufficiently large $k$. (\ref{a-d-k}) is
weaker than (\ref{Lion}). If $u\equiv 0$, we already see that (\ref{Lion}) is weaker than (\ref{Tr}), and obviously
(\ref{a-d-k}) is stronger than (\ref{Tr}).

A natural question is to find the high dimensional analogue of (\ref{a-d}). Let $\Omega$ be a smooth bounded domain in
$\mathbb{R}^n$ $(n\geq 3)$. We proved in \cite{Yang-JFA}
that for any $0\leq \alpha<\lambda_1(\Omega)$,
\be\label{jfa}\sup_{u\in W_0^{1,n}(\Omega),\,\|\nabla u\|_n^n\leq 1}\int_\Omega
e^{\alpha_n|u|^{\f{n}{n-1}}(1+\alpha\|u\|_n^n)^{\f{1}{n-1}}}dx<\infty,\ee
and that the supremum is infinite when $\alpha\geq \lambda_1(\Omega)$,
where $\alpha_n=n\omega_{n-1}^{1/(n-1)}$, $\omega_{n-1}$ is the area of the unit sphere in $\mathbb{R}^n$,
and $\lambda_1(\Omega)$ is defined by
$$\lambda_1(\Omega)=\inf_{u\in W_0^{1,n}(\Omega),\,u\not\equiv 0}\f{\int_\Omega |\nabla u|^ndx}{\int_\Omega |u|^ndx}.$$

Trudinger-Moser inequalities on Riemannian manifolds were due to T. Aubin \cite{Aubin}, J. Moser \cite{14}, P. Cherrier \cite{Cherrier1,Cherrier2},
  and L. Fontana \cite{Fontana}. Also a few results was recently obtained, on  complete noncompact Riemannian manifolds, by
G. Mancini and K. Sandeep \cite{M-S1,M-S2} and the author \cite{Yang-jfa}. One may ask whether or not the analogue of (\ref{a-d}) holds on compact Riemannian surface.
 Let $(\Sigma,g)$ be a compact Riemannian surface without boundary. In \cite{Yang-Tran}, we proved
 the following: For any $\alpha$, $0\leq \alpha<\lambda_1(\Sigma)$, there holds
\be\label{tran}\sup_{u\in W^{1,2}(\Sigma),\,\|\nabla_g u\|_2\leq 1,\,\int_\Sigma udv_g=0}\int_{\Sigma}e^{4\pi u^2(1+\alpha\|u\|_2^2)}
<\infty,\ee
and the supremum is infinite when $\alpha\geq \lambda_1(\Sigma)$, where $W^{1,2}(\Sigma)$ is the usual Sobolev space and
$\lambda_1(\Sigma)$ is defined by
\be\label{eigen-sigma}\lambda_1(\Sigma)=\inf_{u\in W^{1,2}(\Sigma),\,\int_\Sigma udv_g=0,\,u\not\equiv 0}\f{\int_\Sigma|\nabla_gu|^2dv_g}
{\int_\Sigma u^2dv_g}.\ee
If $(\Sigma,g)$ is a compact Riemannian surface with smooth boundary, the trace Trudinger-Moser inequalities were also established in
\cite{Li-Liu,Yang-MZ}.

Existence of extremal functions for Trudinger-Moser inequality (\ref{Tr}) was first obtained by L. Carleson and A. Chang \cite{CC}
when $\Omega$ is a unit ball. This result was extended by M. Struwe \cite{Struwe} to domains close to
a disc in a measure sense,  and by M. Flucher and K. Lin \cite{Flucher,Lin} to general bounded smooth domains.
Later these results were extended by B. Ruf \cite{Ruf} and Li-Ruf \cite{Li-Ruf} to the whole Euclidean space.
The existence result on compact Riemannian manifold was first obtained by Y. Li \cite{Lijpde}, then by Y. Li and P. Liu \cite{Li-Liu},
 and by the author \cite{Yang-IJM}. For existence of extremal functions for Trudinger-Moser inequality of Adimurthi-Druet type
 (Trudinger-Moser inequalities analogous to (\ref{a-d}) above or (\ref{tint}) below), we
 proved in \cite{Yang-Tran,Lu-Yang} that supremums in (\ref{a-d}) and (\ref{tran}) are attained for sufficiently small
 $\alpha\geq 0$, and that the supremum in (\ref{jfa}) $(n\geq 3)$ is attained for all $\alpha$, $0\leq\alpha<\lambda_1(\Omega)$.
 In this direction, M. de Souza and J. M. do \'O \cite{de-dO} generalized (\ref{a-d}) to the whole Euclidean space
 $\mathbb{R}^2$, and the existence of extremal functions was also obtained.

 Recently G. Wang and D. Ye \cite{Wang-Ye} proved the existence of extremal functions for a singular
 Trudinger-Moser inequality.
 Precisely, let $\mathbb{B}$ be a unit disc in $\mathbb{R}^2$, there holds
   \be\label{W-Y}\sup_{\int_{\mathbb{B}}|\nabla u|^2dx-\int_{\mathbb{B}}\f{u^2}{(1-|x|^2)^2}dx\leq 1}\int_\Omega e^{4\pi u^2}dx<\infty,\ee
   and the supremum is attained.
   Another Trudinger-Moser inequality with interior singularity had been established by
 Adimurthi-Sandeep \cite{A-S} on bounded smooth domain and Adimurthi and the author \cite{Adi-Yang} on the whole Euclidean space.
 Moreover C. Tintarev \cite{Tint} modified the classical Trudinger-Moser inequality as follows: Let $\Omega$ be a
 smooth bounded domain in $\mathbb{R}^2$. There holds
 \be\label{tint}\sup_{\int_\Omega |\nabla u|^2dx-\int_\Omega V(x)u^2dx\leq 1}\int_\Omega e^{4\pi u^2}dx<\infty\ee
 for some class of $V(x)>0$ including (\ref{a-d}) and (\ref{W-Y}).
 For extremal functions for Trudinger-Moser inequalities
 on the hyperbolic space, we refer the reader to G. Mancini, K. Sandeep and C. Tintarev \cite{MST}
 and the references therein.\\

 One of our goals in the current paper is to prove that the supremum in (\ref{tint}) is attained in case $V(x)\equiv\alpha$ with $0\leq \alpha<\lambda_1(\Omega)$. Also we consider similar problem for $0\leq\alpha< \lambda_{\ell+1}(\Omega)$,
 the $(\ell+1)$th eigenvalue of the Laplacian operator with respect to Dirichlet boundary condition.
 Moreover the Riemannian surface case are discussed.
 Our method is combining Carleson-Chang's result \cite{CC} with blow-up analysis. For earlier works involving this method, we refer the reader to
 \cite{Li-Liu-Yang,Yang-JFA,Li-Ruf,Lu-Yang,Wang-Ye}. Before ending this section, we remark that for results in this paper,
 there is a possibility of another proof, which is based on the explicit structure of
 putative weakly vanishing maximizing sequences as concentrating Moser functions. For details about this new method, we refer the reader to Adimurthi and
 C. Tintarev \cite{A-T}.

\section{Main Results}

In this paper we concern extremal functions for Trudinger-Moser inequalities of Adimurthi-Druet type.
Let us first consider the Euclidean case. Let $\Omega$ be a smooth bounded domain in $\mathbb{R}^2$ and $\lambda_1(\Omega)$ be the first
eigenvalue of the Laplacian operator with respect to Dirichlet boundary condition. Denote
\be\label{1alpha}\|u\|_{1,\alpha}=\le(\int_\Omega|\nabla u|^2dx-\alpha\int_\Omega u^2dx\ri)^{1/2}\ee
for any $u\in W_0^{1,2}(\Omega)$ with $\int_\Omega|\nabla u|^2dx-\alpha\int_\Omega u^2dx\geq 0$. Clearly $\|\cdot\|_{1,\alpha}$ is equivalent to the Sobolev norm
$\|\cdot\|_{W_0^{1,2}(\Omega)}$ when $0\leq\alpha<\lambda_1(\Omega)$.
Our first result can be stated as follows:
\begin{theorem} Let $\Omega$ be a smooth bounded domain in $\mathbb{R}^2$,
$\lambda_{1}(\Omega)$ be the first eigenvalue of the Laplacian operator with respect to
Dirichlet boundary condition.
If $0\leq\alpha<\lambda_{1}(\Omega)$, then the supremum
\be\label{eucl}\sup_{u\in W_0^{1,2}(\Omega),\,\|u\|_{1,\alpha}\leq 1}\int_\Omega e^{4\pi u^2}dx\ee
can be attained by some function $u_0\in W_0^{1,2}(\Omega)\cap C^1(\overline{\Omega})$ with $\|u_0\|_{1,\alpha}=1$,
where $\|\cdot\|_{1,\alpha}$ is defined as in (\ref{1alpha}).
\end{theorem}
Theorem 1 obviously implies C. Tintarev's inequality (\ref{tint}) in the case $V(x)\equiv \alpha$, and whence
leads to Adimurthi and O. Druet's original inequality (\ref{a-d}). It should be remarked that Theorem 1 does not imply that
the supremum in (\ref{a-d})
is attained for all $\alpha$, $0\leq \alpha<\lambda_1(\Omega)$. Indeed, (\cite{Lu-Yang},
Theorem 1.2, the case $p=2$) has not been improved so far.
When $\alpha=0$, Theorem 1 recovers the results of L. Carleson and A. Chang \cite{CC},
M. Struwe \cite{Struwe}, M. Flucher \cite{Flucher} and K. Lin \cite{Lin}
in dimension two.\\

Obviously the supremum (\ref{eucl}) is infinite if $\alpha\geq\lambda_1(\Omega)$.
It is natural to ask what we can say when other eigenvalues of the Laplacian operator
are involved. Precisely, let $\lambda_1(\Omega)<\lambda_2(\Omega)<\cdots$ be all distinct eigenvalues of the Laplacian
operator
with respect to Dirichlet boundary condition
and $E_{\lambda_j(\Omega)}$'s be associated eigenfunction spaces, namely
$$E_{\lambda_j(\Omega)}=\le\{u\in W_0^{1,2}(\Omega): -\Delta u=\lambda_j(\Omega)u\ri\}.$$
Note that $W_0^{1,2}(\Omega)$ is a Hilbert space when it is equipped with the inner product
$$\langle u,v\rangle=\int_\Omega \nabla u\nabla vdx,\,\,\,\forall u,v\in W_0^{1,2}(\Omega).$$
For any positive integer $\ell$, We set
$$E_\ell=E_{\lambda_1(\Omega)}\oplus E_{\lambda_2(\Omega)}\oplus\cdots \oplus E_{\lambda_\ell(\Omega)}$$
and
\be\label{eper}E_\ell^\perp=\le\{u\in W_0^{1,2}(\Omega): \int_\Omega  u  vdx=0, \forall v\in E_\ell\ri\}.\ee
It is clear that
$$W_0^{1,2}(\Omega)=E_\ell\oplus E_\ell^\perp,\quad\forall \ell=1, 2,\cdots.$$

Similar to Theorem 1, we have the following:
\begin{theorem} Let $\Omega$ be a smooth bounded domain in $\mathbb{R}^2$, $\ell$ be any positive integer,
$\lambda_{\ell+1}(\Omega)$ be the $(\ell+1)$th eigenvalue of the Laplacian operator with respect to
Dirichlet boundary condition, and $E_\ell^\perp$ be a function space defined as in (\ref{eper}).
Then for any $\alpha$, $0\leq\alpha<\lambda_{\ell+1}(\Omega)$, the supremum
\be\label{su-perp}\sup_{u\in E_\ell^\perp,\,\|u\|_{1,\alpha}\leq 1}\int_\Omega e^{4\pi u^2}dx\ee
can be attained by some $u_0\in E_{\ell}^\perp\cap C^1(\overline{\Omega})$ with $\|u_0\|_{1,\alpha}=1$,
where $\|\cdot\|_{1,\alpha}$ is defined as in (\ref{1alpha}).
\end{theorem}

A quite interesting case of Theorem 2 is $\alpha=0$. It follows that for any positive integer $\ell$,
the supremum
$$\sup_{u\in E_\ell^\perp,\,\|\nabla u\|_{2}\leq 1}\int_\Omega e^{4\pi u^2}dx$$
can be attained by some $u_0\in E_\ell^\perp$ with $\|\nabla u_0\|_2= 1$, which is new so far.
If we denote $E_0=\{0\}$ and $E_0^\perp=W_0^{1,2}(\Omega)$, then Theorem 1 is exactly Theorem 2 in case that $\ell=0$.
\\

Now we consider the manifold case. Let $(\Sigma,g)$ be a compact Riemannian surface without boundary, $\nabla_g$ and $\Delta_g$
be its gradient operator and Laplace-Beltrami operator respectively, and
$\lambda_1(\Sigma)$ be the first eigenvalue of $\Delta_g$ (see (\ref{eigen-sigma}) above).
We denote
\be\label{nr-mani}\|u\|_{1,\alpha}=\le(\int_\Sigma |\nabla_gu|^2dv_g-\alpha\int_\Sigma u^2dv_g\ri)^{1/2}\ee
for all $u\in W^{1,2}(\Sigma)$ with $\int_\Sigma |\nabla_gu|^2dv_g-\alpha\int_\Sigma u^2dv_g\geq 0$.
Now we state an analogue of Theorem 1 as follows:
\begin{theorem}
Let $(\Sigma,g)$ be a compact Riemannian surface without boundary, $\lambda_1(\Sigma)$ be
the first eigenvalue of the Laplace-Beltrami operator.
If $0\leq \alpha<\lambda_{1}(\Sigma)$, then the supremum
$$\sup_{u\in W^{1,2}(\Sigma),\,\int_\Sigma udv_g=0,\, \|u\|_{1,\alpha}\leq 1}\int_\Sigma e^{4\pi u^2}dv_g$$
can be attained by some $u_0\in W^{1,2}(\Sigma)\cap C^1(\Sigma)$ with $\int_\Sigma u_0dv_g=0$ and
$\|u_0\|_{1,\alpha}=1$, where $\|\cdot\|_{1,\alpha}$ is defined as in (\ref{nr-mani}).
\end{theorem}

In case $\alpha=0$, Theorem 3 reduces to a result of Y. Li \cite{Lijpde}. Also it should be remarked  that when
$0\leq\alpha<\lambda_1(\Sigma)$, the inequality
\be\label{sup-1}\Lambda_{1,\alpha}=\sup_{u\in W^{1,2}(\Sigma),\,\int_\Sigma udv_g=0,\, \|u\|_{1,\alpha}\leq 1}\int_\Sigma e^{4\pi u^2}dv_g<+\infty\ee
is stronger than that
\be\label{sup-2}\Lambda_\alpha=\sup_{u\in W^{1,2}(\Sigma),\,\int_\Sigma udv_g=0,\, \|\nabla_gu\|_{2}\leq 1}\int_\Sigma e^{4\pi u^2
(1+\alpha\|u\|_2^2)}dv_g<+\infty,\ee
which was studied by the author in \cite{Yang-Tran}.
In fact, if $u\in W^{1,2}(\Sigma)$ satisfies $\int_\Sigma udv_g=0$ and $\|\nabla_gu\|_2\leq 1$,
then $\|u\|_{1,\alpha}^2\leq 1-\alpha\|u\|_2^2$. Since $1+a\leq \f{1}{1-a}$ for all $a<1$, it follows from (\ref{sup-1}) that
$$\int_\Sigma e^{4\pi u^2(1+\alpha\|u\|_2^2)}dv_g\leq \int_\Sigma e^{4\pi \f{u^2}{1-\alpha\|u\|_2^2}}dv_g\leq
\int_\Sigma e^{4\pi \f{u^2}{\|u\|_{1,\alpha}^2}}dv_g\leq \Lambda_{1,\alpha}.$$
Hence we have $\Lambda_\alpha\leq \Lambda_{1,\alpha}$. This was also observed by C. Tintarev \cite{Tint} in the Euclidean case.
But we caution the reader that Theorem 3 does not imply the existence of extremal functions for (\ref{sup-2}). So it is still open whether or not
extremal functions for (\ref{sup-2}) exit for all $0\leq \alpha<\lambda_1(\Sigma)$.\\

Let $\lambda_1(\Sigma)<\lambda_2(\Sigma)<\cdots$ be all distinct eigenvalues of the Laplace-Beltrami operator $\Delta_g$,
and $E_{\lambda_i(\Sigma)}$'s be
associated eigenfunction spaces, namely
$$ E_{\lambda_i(\Sigma)}=\le\{u\in W^{1,2}(\Sigma): \Delta_g u=\lambda_i(\Sigma)u\ri\},\quad i=1,2,\cdots.$$
For any positive integer $\ell$ we write
$$E_\ell=E_{\lambda_1(\Sigma)}\oplus E_{\lambda_2(\Sigma)}\oplus\cdots\oplus E_{\lambda_\ell(\Sigma)}$$
and
\be\label{eper-1}E_\ell^\perp=\le\{u\in W^{1,2}(\Sigma): \int_\Sigma  u  vdv_g=0, \forall v\in E_\ell\ri\}.\ee
Similar to Theorem 2, we have the following:
\begin{theorem}
Let $(\Sigma,g)$ be a compact Riemannian surface without boundary, $\ell$ be any positive integer,
$\lambda_{\ell+1}(\Sigma)$ be the $(\ell+1)$th eigenvalue of the Laplace-Beltrami operator,
and $E_\ell^\perp$ be a function space defined as in (\ref{eper-1}).
Then for any $\alpha$, $0\leq\alpha<\lambda_{\ell+1}(\Sigma)$, the supremum
$$\sup_{u\in E_\ell^\perp,\,\int_\Sigma udv_g=0,\,\|u\|_{1,\alpha}\leq 1}\int_\Sigma e^{4\pi u^2}dv_g$$
can be attained by some $u_0\in E_{\ell}^\perp\cap C^1(\Sigma)$ with $\int_\Sigma u_0dv_g=0$ and $\|u_0\|_{1,\alpha}=1$,
where $\|\cdot\|_{1,\alpha}$ is defined as in (\ref{nr-mani}).
\end{theorem}

It would be also interesting to find extremal functions for improved trace Trudinger-Moser inequality on compact Riemannian
surface with smooth boundary by blow-up analysis. We would not treat this issue here, but refer the reader to
B. Osgood, R. Phillips and P. Sarnak \cite{OPS}, P. Liu \cite{Liu}, Y. Li and P. Liu \cite{Li-Liu}, and the author \cite{Yang-MZ}
for its development.\\

The proofs of Theorems 1 to 4 are all based on a result of Carleson-Chang \cite{CC} and blow-up analysis.
Pioneer works related to this procedure can be found in Ding et al
\cite{DJLW}, Adimurthi and M. Struwe \cite{Adi-Stru}, Y. Li \cite{Lijpde}, Adimurthi and O. Druet \cite{A-D}. Throughout this paper,
$o_j(1)$ denotes the infinitesimal as $j\ra\infty$, $o_\epsilon(1)$ denotes the infinitesimal as $\epsilon\ra 0$, and so on.
In addition we do not distinguish sequence and subsequence, the reader can recognize it easily from the context. Before
 ending this section, we quote Carleson-Chang's result \cite{CC} for our use later: \\

 \noindent{\it Lemma 5} (Carleson-Chang).  {\it Let $\mathbb{B}$ be the unit disc in $\mathbb{R}^2$.
          Assume $\{v_\epsilon\}_{\epsilon>0}$ is a sequence of functions in $W_0^{1,2}(\mathbb{B})$
          with $\int_{\mathbb{B}}|\nabla v_\epsilon|^2dx=1$. If $|\nabla v_\epsilon|^2dx\rightharpoonup
          \delta_0$ as $\epsilon\ra 0$
          weakly in sense of measure. Then
          $\limsup_{\epsilon\ra 0}\int_{\mathbb{B}}(e^{4\pi v_\epsilon^2}-1)dx\leq\pi
          e$.}\\

The remaining part of this paper is organized as follows: In Section 3, we deal with the Euclidean case and prove Theorems 1 and 2;
In Section 4, we deal with the case of manifold without boundary and prove Theorems 3 and 4.

\section{The Euclidean case}

In this section, using Carleson-Chang's result (Lemma 5) and blow-up analysis,
we prove Theorems 1 and 2. Since the procedure is now standard \cite{Lu-Yang} (for earlier works, see
\cite{DJLW,Adi-Stru,Lijpde,A-D}), we give the outline of the proof and emphasize the difference between our
case and the previous ones. In particular, the essential difference between the proofs of Theorem 1 and (\cite{Lu-Yang}, Theorem 1.2)
is the test function computation in the final step. In the proof of Theorem 2, since the maximizers $u_\epsilon$'s may change signs,
hence Gidas-Ni-Nirenberg's result \cite{GNN} can not be applied to our case. However we can exclude the possibility of
boundary blow-up via Agmon's regularity theorem (\cite{Agmon}, page 444) in an indirect way. In the final step (test function computation),
we must ensure that those test functions belong to the space $E_\ell^\perp$, which is different from the counterpart of the proof of Theorem 1.
 \\

{\it Proof of Theorem 1.} Let $\alpha$ be fixed with $0\leq\alpha<\lambda_1(\Omega)$. We divide the proof into several steps
as following:\\

\noindent {\it Step 1. Maximizers for subcritical functionals}\\

 In this step, we shall prove that
 for any $0<\epsilon<4\pi$, there exists some $u_\epsilon\in W_0^{1,2}(\Omega)\cap C^1(\overline{\Omega})$ with
 $\|u_\epsilon\|_{1,\alpha}=1$ such that
 \be\label{subcrit}\int_\Omega e^{(4\pi-\epsilon)u_\epsilon^2}dx=
 \sup_{u\in W_0^{1,2}(\Omega),\,\|u\|_{1,\alpha}\leq 1}\int_\Omega e^{(4\pi-\epsilon)u^2}dx,\ee
 where $\|\cdot\|_{1,\alpha}$ is defined as in (\ref{1alpha}). Here we do not assume in advance the above supremum is
 finite.
 \\

This is based on a direct method in the calculus of variations. For any $0<\epsilon<4\pi$, we take a sequence
of functions $u_j\in W_0^{1,2}(\Omega)$ verifying that
\be\label{1}\int_\Omega |\nabla u_j|^2dx-\alpha
\int_\Omega u_j^2dx\leq 1\ee
and that as $j\ra\infty$,
\be\label{2}
\int_\Omega e^{(4\pi-\epsilon)u_j^2}dx\ra \sup_{u\in W_0^{1,2}(\Omega),\,\|u\|_{1,\alpha}\leq 1}\int_\Omega e^{(4\pi-\epsilon)u^2}dx.
\ee
 It follows from (\ref{1}) and $0\leq\alpha<\lambda_1(\Omega)$ that $u_j$ is bounded in $W_0^{1,2}(\Omega)$.
Thus we can assume up to a subsequence, $u_j\rightharpoonup u_\epsilon$  weakly  in
  $W_0^{1,2}(\Omega)$,  $u_j\ra u_\epsilon$  strongly in  $L^2(\Omega)$, and
  $u_j\ra u_\epsilon$ a.e.  in $\Omega$.
  Clearly we have that
  \be\label{ueps}0\leq\int_\Omega |\nabla u_\epsilon|^2dx-\alpha\int_\Omega u_\epsilon^2dx\leq \liminf_{j\ra\infty}
  \le(\int_\Omega|\nabla u_j|^2dx-\alpha\int_\Omega u_j^2dx\ri)\leq 1\ee
  and that
  \bea\nonumber\int_\Omega|\nabla u_j-\nabla u_\epsilon|^2dx&=& \int_\Omega|\nabla u_j|^2dx-\int_\Omega|\nabla u_\epsilon|^2dx
  +o_j(1)\\\label{less1}
  &\leq&1-\int_\Omega|\nabla u_\epsilon|^2dx+\alpha\int_\Omega u_\epsilon^2dx+o_j(1).\eea
   Combining (\ref{ueps}) and (\ref{less1}), we conclude
  $$\limsup_{j\ra\infty}\int_\Omega|\nabla u_j-\nabla u_\epsilon|^2dx\leq 1.$$
  It follows
  from Lion's inequality (\ref{Lion}) that
  $e^{(4\pi-\epsilon)u_j^2}$ is bounded in $L^q(\Omega)$ for some $q>1$. Hence $e^{(4\pi-\epsilon)u_j^2}\ra e^{(4\pi-\epsilon)u_\epsilon^2}$
  strongly in $L^1(\Omega)$. This together with (\ref{2}) immediately leads to (\ref{subcrit}). Obviously the supremum in (\ref{subcrit})
  is strictly greater than $|\Omega|$, the volume of $\Omega$. Therefore $u_\epsilon\not\equiv 0$. If $\|u_\epsilon\|_{1,\alpha}<1$,
  we set $\widetilde{u}_\epsilon=
  u_\epsilon/\|u_\epsilon\|_{1,\alpha}$, then we obtain $\|\widetilde{u}_\epsilon\|_{1,\alpha}=1$ and
  $$\sup_{u\in W_0^{1,2}(\Omega),\,\|u\|_{1,\alpha}\leq 1}\int_\Omega e^{(4\pi-\epsilon)u^2}dx\geq
  \int_\Omega e^{(4\pi-\epsilon)\widetilde{u}_\epsilon^2}dx>\int_\Omega e^{(4\pi-\epsilon){u}_\epsilon^2}dx.$$
  This contradicts (\ref{subcrit}). Hence $\|{u}_\epsilon\|_{1,\alpha}=1$.

It is not difficult to see that $u_\epsilon$ satisfies the Euler-Lagrange equation
\be\label{E-L}\le\{
  \begin{array}{lll}
  -\Delta u_\epsilon-\alpha u_\epsilon=\f{1}{\lambda_\epsilon}u_\epsilon e^{(4\pi-\epsilon)u_\epsilon^2}\,\,\,{\rm in}\,\,\,
  \Omega,\\[1.5ex] u_\epsilon>0\,\,\,{\rm in}\,\,\,
  \Omega,\\[1.5ex] \int_\Omega|\nabla u_\epsilon|^2dx-\alpha\int_\Omega u_\epsilon^2dx=1,\\[1.5ex]
  \lambda_\epsilon=\int_\Omega u_\epsilon^2 e^{(4\pi-\epsilon)u_\epsilon^2}dx.
  \end{array}
  \ri.\ee
  Applying elliptic estimates to (\ref{E-L}), we have $u_\epsilon\in C^1(\overline{\Omega})$.
  Let $c_\epsilon=u_\epsilon(x_\epsilon)=\max_\Omega u_\epsilon$. Since $\|u_\epsilon\|_{1,\alpha}=1$,
  without loss of generality, we assume $u_\epsilon$ converges to $u^\ast$  weakly in $W_0^{1,2}(\Omega)$,
  strongly in $L^2(\Omega)$, and almost everywhere in $\Omega$.
  If $c_\epsilon$ is bounded, then $e^{(4\pi-\epsilon)u_\epsilon^2}$ is bounded in $L^\infty(\Omega)$, and thus
  $e^{(4\pi-\epsilon)u_\epsilon^2}$ converges to $e^{4\pi {u^\ast}^2}$ in $L^1(\Omega)$. Hence for any $u\in W_0^{1,2}(\Omega)$
  with $\|u\|_{1,\alpha}\leq 1$, we have by (\ref{subcrit}) that
  $$\label{c-bounded}\int_\Omega e^{4\pi u^2}dx=\lim_{\epsilon\ra 0}\int_\Omega e^{(4\pi-\epsilon) u^2}dx
  \leq \lim_{\epsilon\ra 0}\int_\Omega e^{(4\pi-\epsilon) u_\epsilon^2}dx=\int_\Omega e^{4\pi {u^\ast}^2}dx.$$
  This implies that
  \be\label{bd-extre}\sup_{u\in W_0^{1,2}(\Omega),\,\|u\|_{1,\alpha}\leq 1}\int_\Omega e^{4\pi u^2}dx= \int_\Omega e^{4\pi
  {u^\ast}^2}dx.\ee
  So $u^\ast\in W_0^{1,2}(\Omega)$ attains the above supremum. Obviously $\|u^\ast\|_{1,\alpha}=1$. Applying elliptic estimates to its
  Euler-Lagrange equation, we obtain $u^\ast\in C^1(\overline\Omega)$. Therefore $u^\ast$ is the desired extremal function.
    Hence we assume $c_\epsilon\ra\infty$ in the sequel.
  Without loss of generality, we assume $x_\epsilon\ra x_0\in\overline\Omega$.
  By a result of Gidas-Ni-Nirenberg (\cite{GNN}, page 223), the distance between $x_\epsilon$ and $\p \Omega$ must be greater than $\delta>0$
  depending only on $\Omega$. Therefore $x_0\not\in\p\Omega$.
  \\


\noindent{\it Step 2. Energy concentration phenomenon}\\

{ In this step we shall prove that $u_\epsilon\rightharpoonup 0$ weakly in $W_0^{1,2}(\Omega)$, $u_\epsilon \ra 0$
strongly in $L^q(\Omega)$ for any $q> 1$, and $|\nabla u_\epsilon|^2dx\rightharpoonup\delta_{x_0}$ weakly in sense of measure
as $\epsilon\ra 0$,
where $\delta_{x_0}$ is the usual Dirac measure centered at $x_0$.}\\

Noting that $\|u_\epsilon\|_{1,\alpha}=1$, we can assume $u_\epsilon\rightharpoonup u_0$ weakly in $W_0^{1,2}(\Omega)$,
and $u_\epsilon\ra u_0$ strongly in $L^q(\Omega)$ for any $q>1$. It follows that
\be\label{3}\int_\Omega|\nabla u_\epsilon|^2dx=1+\alpha\int_\Omega u_0^2dx+o(1),\ee
and that
\be\label{4}\int_\Omega|\nabla (u_\epsilon-u_0)|^2dx=1-\int_\Omega|\nabla u_0|^2dx+\alpha\int_\Omega u_0^2dx+o(1).\ee
Suppose $u_0\not\equiv0$. In view of (\ref{4}), Lions' inequality (\ref{Lion}) implies that $e^{4\pi u_\epsilon^2}$ is bounded in $L^q(\Omega)$
for any fixed $q$ with $1<q<1/(1-\|u_0\|_{1,\alpha}^2)$. Then applying
elliptic estimates to (\ref{E-L}), we have that $u_\epsilon$ is uniformly bounded in $\Omega$, which contradicts $c_\epsilon\ra\infty$.
Therefore $u_0\equiv0$ and (\ref{3}) becomes
\be\label{5}\int_\Omega|\nabla u_\epsilon|^2dx=1+o_\epsilon(1).\ee
Suppose $|\nabla u_\epsilon|^2dx\rightharpoonup \mu$ in sense of measure. If $\mu\not =\delta_{x_0}$, then in view of (\ref{5}) and $u_0\equiv0$,
 we can choose some $r_0>0$
and a cut-off
function $\phi\in C_0^1(\mathbb{B}_{r_0}(x_0))$, which is equal to $1$ on $\mathbb{B}_{r_0/2}(x_0)$, such that
$\mathbb{B}_{r_0}(x_0)\subset\Omega$ and
$$\limsup_{\epsilon\ra 0}\int_{\mathbb{B}_{r_0}(x_0)}|\nabla(\phi u_\epsilon)|^2dx<1.$$
 By the classical Trudinger-Moser inequality (\ref{Tr}), $e^{(4\pi-\epsilon)(\phi u_\epsilon)^2}$ is bounded in
 $L^r(\mathbb{B}_{r_0}(x_0))$ for some $r>1$.
 Applying elliptic estimates to (\ref{E-L}), we have that $u_\epsilon$ is uniformly bounded in $\mathbb{B}_{r_0/2}(x_0)$,
 which contradicts $c_\epsilon\ra\infty$ again.
 Therefore $|\nabla u_\epsilon|^2dx\rightharpoonup \delta_{x_0}$ and Step 2 is finished.\\

\noindent{\it Step 3. Blow-up analysis for $u_\epsilon$}\\

We set
$$\label{r-eps} r_\epsilon=\sqrt{\lambda_\epsilon}c_\epsilon^{-1}e^{-(2\pi-\epsilon/2)c_\epsilon^2}.$$
 For any $0<\delta<4\pi$, we have by using the H\"older inequality and the classical Trudinger-Moser inequality (\ref{Tr}),
 $$\label{l-eps}\lambda_\epsilon=\int_\Omega u_\epsilon^2 e^{(4\pi-\epsilon)u_\epsilon^2}dx\leq e^{\delta c_\epsilon^2}
 \int_\Omega u_\epsilon^2 e^{(4\pi-\epsilon-\delta)u_\epsilon^2}dx\leq Ce^{\delta c_\epsilon^2}$$
 for some constant $C$ depending only on $\delta$.
 This leads to
 \be\label{r-0} r_\epsilon^2\leq C c_\epsilon^{-2}e^{-(4\pi-\epsilon-\delta)c_\epsilon^2}\ra 0\quad{\rm as}\quad \epsilon\ra 0.\ee
 Let
 $$\Omega_\epsilon=\{x\in \mathbb{R}^2:x_\epsilon+r_\epsilon x\in\Omega\}.$$
   Define two blow-up sequences of functions on $\Omega_\epsilon$ as
 $$
   \psi_\epsilon(x)=c_\epsilon^{-1}u_\epsilon(x_\epsilon+r_\epsilon
   x),\quad
   \varphi_\epsilon(x)=c_\epsilon(u_\epsilon(x_\epsilon+r_\epsilon
   x)-c_\epsilon).
   $$
   A direct computation shows
   \be\label{p-s-eq}
   -\la \psi_\epsilon=\alpha r_\epsilon^2\psi_\epsilon+c_\epsilon^{-2}\psi_\epsilon e^{(4\pi-\epsilon)(u_\epsilon^2-c_\epsilon^2)}\quad{\rm
   in}\quad \Omega_\epsilon,
   \ee
   \be\label{phi-eq}
   -\la\varphi_\epsilon=\alpha r_\epsilon^2c_\epsilon^2\psi_\epsilon+\psi_\epsilon e^{(4\pi-\epsilon)(1+\psi_\epsilon)\varphi_\epsilon}\quad{\rm
   in}\quad \Omega_\epsilon.
   \ee
   We now investigate the convergence behavior of $\psi_\epsilon$ and $\varphi_\epsilon$.
   Note that $\Omega_\epsilon\ra \mathbb{R}^2$ as $\epsilon\ra 0$. Since $|\psi_\epsilon|\leq 1$ and
   $\Delta\psi_\epsilon(x)\ra 0$ uniformly in $x\in\Omega_\epsilon$ as $\epsilon\ra 0$, we have by elliptic estimates that
   $\psi_\epsilon\ra \psi$ in $C^1_{\rm loc}(\mathbb{R}^2)$, where $\psi$ is a bounded harmonic function in $\mathbb{R}^2$.
   Note that $\psi(0)=\lim_{\epsilon\ra 0}\psi_\epsilon(0)=1$. The Liouville theorem implies that $\psi\equiv 1$ on
   $\mathbb{R}^2$. Thus we have
   \be\label{ps-t-1}\psi_\epsilon\ra 1\quad{\rm in}\quad C^1_{\rm loc}(\mathbb{R}^2).\ee
   By (\ref{r-0}), we have $r_\epsilon^2c_\epsilon^2\ra 0$ as $\epsilon\ra 0$. Note also that
   $$\varphi_\epsilon(x)\leq \varphi_\epsilon(0)=0\,\,\,{\rm for\,\,all}\,\,\,
   x\in \Omega_\epsilon.$$
   Thus $\Delta \varphi_\epsilon$ is uniformly bounded in $\Omega_\epsilon$. We then conclude by applying elliptic estimates to
   the equation (\ref{phi-eq}) that
   \be\label{phibubble}
   \varphi_\epsilon\ra \varphi\quad{\rm in}\quad
   C^1_{\rm{loc}}(\mathbb{R}^2),
   \ee
   where $\varphi$ satisfies
   $$\label{bubbel}
   \le\{\begin{array}{lll}&\la
     \varphi=-e^{8\pi\varphi}\quad{\rm in}\quad\mathbb{R}^2\\[1.2 ex]
     &\varphi(0)=0=\sup_{\mathbb{R}^2}\varphi\\[1.2 ex] &\int_{\mathbb{R}^2}e^{8\pi\varphi}dx\leq 1.
     \end{array}
    \ri.$$
    By a result of Chen-Li \cite{CL}, we have
    \be\label{v-r}\varphi(x)=-\f{1}{4\pi}\log(1+\pi|x|^2)\ee
    and
    $$\int_{\mathbb{R}^2}e^{8\pi\varphi}dx=1.$$
    To understand the convergence behavior away from the blow-up point $x_0$, we need to investigate how $c_\epsilon u_\epsilon$
     converges. By a repetitive argument of (\cite{Lu-Yang}, Lemma 3.6), we have that
     \be\label{cu}c_\epsilon u_\epsilon\rightharpoonup G\,\,\,{\rm weakly\,\, in}\,\,\, W_0^{1,q}(\Omega),\,\,\,\forall
     1<q<2,\ee
     where $G\in C^1(\overline\Omega\setminus\{x_0\})$
     is the Green function satisfying the equation
     \be\label{Green}
     \left\{\begin{array}{lll}
            -\Delta G-\alpha G=\delta_{x_0}\quad\rm{in}\quad
            \Omega\\[1.5ex]
            G=0\quad {\rm on}\quad \p\Omega.
         \end{array}\right.\ee
         Moreover,
         \be\label{G-C1-1}c_\epsilon u_\epsilon\ra G\quad{\rm in}\quad
         C^1_{\rm loc}(\overline{\Omega}\setminus\{x_0\}).\ee

         \noindent{\it Step 4. Upper bound estimate}\\

          In view of (\ref{Green}) and (\ref{G-C1-1}), $G$ can be represented by
          \be\label{Ax0}
          G=-\f{1}{2\pi}\log
          |x-x_0|+A_{x_0}+\psi_\alpha(x),\ee where $A_{x_0}$ is a
          constant depending on $x_0$ and $\alpha$, $\psi_\alpha\in
          C^1(\Omega)$ and $\psi_\alpha(x_0)=0$.
          This leads to
          \bna
          \int_{\Omega\setminus B_\delta(x_0)}|\nabla
          G|^2dx&=&\alpha\int_{\Omega\setminus
          B_\delta(x_0)}G^2dx+\int_{\p(\Omega\setminus
          B_\delta(x_0))}G\f{\p G}{\p n}ds\\
          &=&\f{1}{2\pi}\log\f{1}{\delta}+A_{x_0}+\alpha\|G\|_2^2+o_\delta(1).
          \ena
          Hence we obtain
          \be\label{Omega-B}
          \int_{\Omega\setminus B_\delta(x_0)}|\nabla
          u_\epsilon|^2dx=\f{1}{c_\epsilon^2}
          \le(\f{1}{2\pi}\log\f{1}{\delta}+A_{x_0}+\alpha\|G\|_2^2+o_\delta(1)+o_\epsilon(1)\ri).
          \ee
          Let $s_\epsilon=\sup_{\p B_\delta(x_0)}u_\epsilon$ and
          $\overline{u}_\epsilon=(u_\epsilon-s_\epsilon)^+$. Then
          $\overline{u}_\epsilon\in W_0^{1,2}(B_\delta(x_0))$.
          By (\ref{Omega-B}) and the fact that $\int_{B_\delta(x_0)}|\nabla
          u_\epsilon|^2dx=1-\int_{\Omega\setminus B_\delta(x_0)}|\nabla
          u_\epsilon|^2dx+\alpha\int_\Omega u_\epsilon^2$, we have
          \be\label{4.4}
          \int_{B_\delta(x_0)}|\nabla
          \overline{u}_\epsilon|^2dx\leq \tau_\epsilon=1-\f{1}{c_\epsilon^2}
          \le(\f{1}{2\pi}\log\f{1}{\delta}+A_{x_0}+o_\delta(1)+o_\epsilon(1)\ri).
          \ee
          This together with Lemma 5 (see the end of Section 2) leads to
          \be\label{B-delta}
          \limsup_{\epsilon\ra 0}\int_{B_\delta(x_0)}(e^{4\pi
          \overline{u}_\epsilon^2/\tau_\epsilon}-1)dx\leq
          \pi\delta^2e.
          \ee
          By (\ref{phibubble}), we have on $B_{Rr_\epsilon}(x_\epsilon)$ that
          $u_\epsilon(x)=c_\epsilon+\f{1}{c_\epsilon}\varphi(\f{x-x_\epsilon}{r_\epsilon})$,
          which together with the fact that $c_\epsilon u_\epsilon\ra
          G$ in $L^2(\Omega)$, gives on
          $B_{Rr_\epsilon}(x_\epsilon)$,
          \bna
          (4\pi-\epsilon)
          u_\epsilon^2&\leq&4\pi(\overline{u}_\epsilon+s_\epsilon)^2\\
          &\leq&
          4\pi\overline{u}_\epsilon^2+8\pi
          s_\epsilon\overline{u}_\epsilon+o_\epsilon(1)\\
          &\leq&4\pi\overline{u}_\epsilon^2-4\log\delta+8\pi
          A_{x_0}+o_\epsilon(1)+o_\delta(1)\\
          &\leq&4\pi\overline{u}_\epsilon^2/\tau_\epsilon-2\log\delta+4\pi
          A_{x_0}+o(1).
          \ena
          Therefore
          \bna
          \int_{B_{Rr_\epsilon}(x_\epsilon)}e^{(4\pi-\epsilon)
          u_\epsilon^2}dx&\leq& \delta^{-2}e^{4\pi
          A_{x_0}+o(1)}\int_{B_{Rr_\epsilon}(x_\epsilon)}e^{4\pi\overline{u}_\epsilon^2/\tau_\epsilon}dx\\
          &=&\delta^{-2}e^{4\pi
          A_{x_0}+o(1)}\int_{B_{Rr_\epsilon}(x_\epsilon)}(e^{4\pi\overline{u}_\epsilon^2/\tau_\epsilon}-1)dx+o(1)\\
          &\leq&\delta^{-2}e^{4\pi
          A_{x_0}+o(1)}\int_{B_\delta(x_0)}(e^{4\pi\overline{u}_\epsilon^2/\tau_\epsilon}-1)dx.
          \ena
          This together with (\ref{B-delta}) leads to
          \be\label{u-p}
          \limsup_{\epsilon\ra 0}\int_{B_{Rr_\epsilon}(x_\epsilon)}e^{(4\pi-\epsilon)
          u_\epsilon^2}dx\leq \pi e^{1+4\pi A_{x_0}}.
          \ee
          By the same argument as in the proof of (\cite{Lu-Yang}, Lemma 3.3), we get
          \be\label{u-p-1}\lim_{\epsilon\ra 0}\int_\Omega e^{(4\pi-\epsilon)u_\epsilon^2}dx\leq |\Omega|+
     \lim_{R\ra+\infty}\limsup_{\epsilon\ra 0}\int_{B_{Rr_\epsilon}(x_\epsilon)}e^{(4\pi-\epsilon)u_\epsilon^2}dx\ee
          Combining (\ref{u-p}) and (\ref{u-p-1}), we conclude
          \be\label{upperbound}
          \sup_{u\in W_0^{1,2}(\Omega),\,\|u\|_{1,\alpha}\leq 1}\int_\Omega e^{4\pi u^2}dx=\limsup_{\epsilon\ra 0}\int_{\Omega}e^{(4\pi-\epsilon)
          u_\epsilon^2}dx\leq |\Omega|+\pi e^{1+4\pi A_{x_0}}.
          \ee

         \noindent {\it Step 5. Existence of extremal functions}\\

   We will construct a sequence of functions
     $\phi_\epsilon\in W_0^{1,2}(\Omega)$ such that $\|\phi_\epsilon\|_{1,\alpha}=1$
     and
     \be\label{blowsequence}
     \int_\Omega e^{4\pi\phi_\epsilon^2}dx>
     |\Omega|+\pi e^{1+4\pi A_{x_0}}
     \ee
     for sufficiently small
     $\epsilon>0$. The contradiction between  (\ref{upperbound}) and
    (\ref{blowsequence}) implies that $c_\epsilon$ must be bounded. Then applying elliptic
     estimates to (\ref{E-L}), we conclude the existence of extremal function and finish the proof of Theorem 1.

     To prove (\ref{blowsequence}),
     we recall (\ref{Ax0}) and write
     $r(x)=|x-x_0|$.     Set
     \be\label{ppp}\phi_\epsilon=\le\{
     \begin{array}{llll}
     &c+\f{-\f{1}{4\pi}\log(1+\pi\f{r^2}{\epsilon^2})+B}{c}
     \quad &{\rm for} &r\leq R\epsilon\\[1.5ex]
     &\f{G-\eta \psi_\alpha}{c}\quad &{\rm for} & R\epsilon<
     r<2R\epsilon\\[1.2ex]
     &\f{G}{c}\quad &{\rm for} & r\geq 2R\epsilon,
     \end{array}
     \ri.\ee
     where $R=-\log\epsilon$, $\eta\in C_0^\infty(B_{2R\epsilon}(x_0))$ verifying that $\eta=1$ on $B_{R\epsilon}(x_0)$ and
     $\|\nabla \eta\|_{L^\infty}
     =O(\f{1}{R\epsilon})$, $B$ is a constant to be determined
     later, and $c$ depending only on $\epsilon$ will also be
     chosen later such that $R\epsilon\ra 0$ and $R\ra +\infty$.
     In order to assure that $\phi_\epsilon\in W_0^{1,2}(\Omega),$ we set
     $$
     c+\f{1}{c}\le(-\f{1}{4\pi}\log(1+\pi R^2)+B\ri)
     =\f{1}{c}\le(-\f{1}{2\pi}\log (R\epsilon)+A_{x_0}\ri),
     $$
     which gives
     \be\label{2pic2-1}
     2\pi c^2=-\log\epsilon-2\pi B+2\pi A_{x_0}+\f{1}{2}\log \pi
     +O(\f{1}{R^2}).
     \ee
     A delicate but straightforward calculation shows
     \bna
     \int_{\Omega}|\nabla \phi_\epsilon|^2dx&=&\f{1}{4\pi c^2}\le(
     2\log\f{1}{\epsilon}+\log\pi-1+4\pi A_{x_0}+4\pi\alpha\|G\|_2^2\ri.\\
     &&\le.+O(\f{1}{R^2})+
     O(R\epsilon\log(R\epsilon))\ri)
     \ena
     and
     $$
     \int_\Omega\phi_\epsilon^2dx=\f{1}{c^2}\le(\int_\Omega G^2dx+O(R\epsilon\log(R\epsilon))\ri),
     $$
     which yields
     \bna \|\phi_\epsilon\|_{1,\alpha}^2&=&\int_\Omega(|\nabla\phi_\epsilon|^2-\alpha\phi_\epsilon^2)dx\\
     &=&
     \f{1}{4\pi c^2}\le(
     2\log\f{1}{\epsilon}+\log\pi-1+4\pi A_{x_0}
     +O(\f{1}{R^2})+
     O(R\epsilon\log(R\epsilon))\ri)\ena
     Set $\|\phi_\epsilon\|_{1,\alpha}=1$, we have
     \be\label{c2-1}
     c^2=-\f{\log\epsilon}{2\pi}+\f{\log\pi}{4\pi}-\f{1}{4\pi}+A_{x_0}
     +O(\f{1}{R^2})+O(R\epsilon\log(R\epsilon)).
     \ee
     It follows from (\ref{2pic2-1}) and (\ref{c2-1}) that
     \be{\label{B-1}}
     B=\f{1}{4\pi}+O(\f{1}{R^2})+O(R\epsilon\log(R\epsilon)).
     \ee
     Clearly we have on $B_{R\epsilon}(x_0)$
     $$4\pi\phi_\epsilon^2\geq 4\pi
     c^2-2\log(1+\pi\f{r^2}{\epsilon^2})+8\pi B.$$
     This together with (\ref{c2-1}) and (\ref{B-1}) yields
     \be\label{BRE-1}
     \int_{B_{R\epsilon}(x_0)} e^{4\pi\phi_\epsilon^2}dx\geq
     \pi e^{1+4\pi A_{x_0}}
     +O(\f{1}{R^2}).
     \ee
     On the other hand,
     \bea
      \label{O-BRE-1}\int_{\Omega\setminus
      B_{R\epsilon}(x_0)}e^{4\pi\phi_\epsilon^2}dx&\geq&\int_{\Omega\setminus
      B_{2R\epsilon}(x_0)}(1+4\pi\phi_\epsilon^2)dx\\{\nonumber}
      &\geq& |\Omega|+4\pi\f{\|G\|_2^2}{c^2}+o(\f{1}{c^2}).
     \eea
     Recalling (\ref{c2-1}) and the choice of $R=-\log\epsilon$,
     we conclude (\ref{blowsequence}) for sufficiently small $\epsilon>0$ by combining (\ref{BRE-1}) and (\ref{O-BRE-1}).
     $\hfill\Box$\\

     Before proving Theorem 2, we state a special version of a regularity theorem due to S. Agmon (\cite{Agmon}, page 444),
     which is essential for excluding boundary blow-up.\\

     \noindent{\it Lemma 6. Let $\Omega$ be a smooth bounded domain in $\mathbb{R}^2$,  $u\in L^r(\Omega)$ for some $r>1$,
     and $f\in L^q(\Omega)$ for some $q>1$. Suppose that for all functions $v\in C^2(\overline{\Omega})\cap W_0^{1,q}(\Omega)$,
     $$\int_\Omega u\Delta vdx=\int_\Omega fvdx.$$
     Then $u\in W^{2,q}(\Omega)\cap W_0^{1,q}(\Omega)$.}\\

     {\it Proof of Theorem 2.}
     Firstly, we fix several notations concerning the function space $E_\ell^\perp$ defined as in (\ref{eper}).
     Let $\lambda_1(\Omega)<\lambda_2(\Omega)<\cdots$ be all distinct eigenvalues of
     the Laplacian with respect to Dirichlet boundary condition, and $E_{\lambda_i(\Omega)}$'s be
     associated eigenfunction spaces. It is known that $\lambda_i(\Omega)\ra +\infty$ as $i\ra+\infty$
     and each space $E_{\lambda_i(\Omega)}$ has finite dimension (see \cite{Brezis}, Theorem 9.31).
     We can assume
     $${\rm dim}E_{\lambda_i(\Omega)}=n_i,\quad i=1,2,\cdots.$$
     Moreover we can find a basis $(e_{ij})$ $(1\leq j\leq n_i,1\leq i\leq\ell)$ of $E_\ell$ verifying
     \be\label{50}\le\{\begin{array}{lll}
     E_{\lambda_i(\Omega)}={\rm span}\{e_{i1},\cdots,e_{in_i}\},\quad i=1,\cdots,\ell,\\
     [1.2ex]E_\ell={\rm span}\{e_{11},\cdots,e_{1n_1}, e_{21},\cdots,e_{2n_2},\cdots,e_{\ell 1},\cdots,e_{\ell n_\ell}\},\\
     [1.2ex] \int_\Omega|e_{ij}|^2dx=1,\\
     [1.2ex] \int_\Omega e_{ij}e_{kl}dx=0,\,\, i\not= k\,\,{\rm or}\,\, j\not=l.
     \end{array}\ri.\ee
     Note that
     $$E_\ell^\perp=\le\{u\in W_0^{1,2}(\Omega): \int_\Omega  u  e_{ij}dx=0, 1\leq j\leq n_i, 1\leq i\leq \ell\ri\}.$$

      Secondly, let $0\leq\alpha<\lambda_{\ell+1}(\Omega)$ be fixed, we shall find maximizers for subcritical Trudinger-Moser functionals.
      Analogous to Step 1 of the proof of Theorem 1,
     for any $\epsilon$, $0<\epsilon<4\pi$, there exists
     some $u_\epsilon\in E_\ell^\perp\cap C^1(\overline{\Omega})$ with
 $\|u_\epsilon\|_{1,\alpha}=1$ such that
 \be\label{subcrit1-1}\int_\Omega e^{(4\pi-\epsilon)u_\epsilon^2}dx=
 \sup_{u\in E_\ell^\perp,\,\|u\|_{1,\alpha}\leq 1}\int_\Omega e^{(4\pi-\epsilon)u^2}dx,\ee
 where $\|\cdot\|_{1,\alpha}$ is defined as in (\ref{1alpha}). Moreover $u_\epsilon$ satisfies the
 Euler-Lagrange equation
  \be\label{E-L-1}\le\{
  \begin{array}{lll}
  -\Delta u_\epsilon-\alpha u_\epsilon=\f{1}{\lambda_\epsilon}u_\epsilon e^{(4\pi-\epsilon)u_\epsilon^2}-
  \sum_{i=1}^\ell\sum_{j=1}^{n_i}\f{\gamma_{ij,\epsilon}}{\lambda_\epsilon}e_{ij}\,\,\,{\rm in}\,\,\,
  \Omega,\\[1.5ex] u_\epsilon\in E_\ell^\perp\cap C^1(\overline{\Omega}),\\[1.5ex]
  \lambda_\epsilon=\int_\Omega u_\epsilon^2 e^{(4\pi-\epsilon)u_\epsilon^2}dx,\\
  [1.5ex] \gamma_{ij,\epsilon}=\int_\Omega e_{ij}u_\epsilon e^{(4\pi-\epsilon)u_\epsilon^2}dx.
  \end{array}
  \ri.\ee
  Without loss of generality we can assume
  \bea
  &&u_\epsilon\rightharpoonup u_0\quad{\rm weakly\,\,\,in}\quad W_0^{1,2}(\Omega),\label{w-c}\\
  [1.2ex] &&u_\epsilon\ra u_0\quad{\rm strongly\,\,\,in}\quad L^p(\Omega), \,\,\,
  \forall p>1,{\label{s-c}}\\
  [1.2ex] &&u_\epsilon\ra u_0\quad{\rm a.\,e.\,\,\,in}\quad \Omega.\label{ae-c}
  \eea
  Since $u_\epsilon\in E_\ell^\perp$, we have by (\ref{s-c})
  $$\int_\Omega u_0e_{ij}dx=\lim_{\epsilon\ra 0}\int_\Omega u_\epsilon e_{ij}dx=0,\,\,\, 1\leq j\leq n_i,
  \,\,1\leq i\leq \ell,$$
  which together with (\ref{w-c}) implies that $u_0\in E_\ell^\perp$ and $\|u_0\|_{1,\alpha}\leq 1$.

  If $u_\epsilon$ is bounded in $C^0(\overline\Omega)$, then for any $v\in E_\ell^\perp$ with
  $\|v\|_{1,\alpha}\leq 1$, (\ref{subcrit1-1}), (\ref{ae-c}) and Lebesgue's dominated convergence theorem
  lead to
  $$\int_\Omega e^{4\pi {v}^2}dx=\lim_{\epsilon\ra 0}\int_\Omega e^{(4\pi-\epsilon){v}^2}dx
  \leq\lim_{\epsilon\ra 0}\int_{\Omega}e^{(4\pi-\epsilon) u_\epsilon^2}dx=\int_\Omega e^{4\pi u_0^2}dx.$$
  Hence we have
  \be\label{b-d}\int_\Omega e^{4\pi u_0^2}dx=\sup_{u\in E_\ell^\perp,\,\|u\|_{1,\alpha}\leq 1}\int_\Omega e^{4\pi u^2}dx.\ee
  It is easy to see that $\|u_0\|_{1,\alpha}=1$.
  Applying elliptic estimates to the Euler-Lagrange equation of $u_0$, we have $u_0\in C^1(\overline{\Omega})$.
  Thus $u_0$ is the desired extremal function.

  In the sequel we assume up to a subsequence
  $$\|u_\epsilon\|_{C^0(\overline\Omega)}=\max_{\overline\Omega}|u_\epsilon|\ra +\infty\quad{\rm as}\quad \epsilon\ra 0.$$

  Thirdly, we perform blow-up analysis. Denote $c_\epsilon=|u_\epsilon(x_\epsilon)|=\|u_\epsilon\|_{C^0(\overline\Omega)}$. Then $c_\epsilon\ra +\infty$ as
  $\epsilon\ra 0$. Without loss of generality
  we assume $c_\epsilon=u_\epsilon(x_\epsilon)$. For otherwise $u_\epsilon$ can be replaced by $-u_\epsilon$ in the
  following blow-up analysis. Then up to a subsequence, $x_\epsilon\ra x_0\in\overline\Omega$. As in Step 2 of the proof
  of Theorem 1, we have $u_0\equiv 0$ and $|\nabla u_\epsilon|^2dx\rightharpoonup \delta_{x_0}$ weakly in sense of measure.
  The only difference is that $\phi u_\epsilon\in W_0^{1,2}(\mathbb{B}_r(x_0)\cap \Omega)$ in case $x_0\in \p\Omega$.\\



  Set
  $$\psi_\epsilon(x)=c_\epsilon^{-1}u_\epsilon(x_\epsilon+r_\epsilon
   x),\quad\varphi_\epsilon(x)=c_\epsilon(u_\epsilon(x_\epsilon+r_\epsilon
   x)-c_\epsilon),\quad x\in\Omega_\epsilon,$$
   where
\be\label{r-eps-1} r_\epsilon=\sqrt{\lambda_\epsilon}c_\epsilon^{-1}e^{-(2\pi-\epsilon/2)c_\epsilon^2}\ee
and
$$\Omega_\epsilon=\{x\in \mathbb{R}^2:x_\epsilon+r_\epsilon x\in\Omega\}.$$
By (\ref{r-0}), we have $r_\epsilon\ra 0$ as $\epsilon\ra 0$. Moreover we claim that up to a subsequence
\be\label{bb}r_\epsilon/{\rm dist}(x_\epsilon,\p\Omega)\ra 0\quad {\rm as}\quad \epsilon\ra 0.\ee
Let $\mathbb{B}$ be a unit disc centered
at $0\in\mathbb{R}^2$. Since $\Omega$ is smooth, we have a a neighborhood $U\subset \mathbb{R}^2$ of $x_0$ and
a bijective map $H:\mathbb{B}\ra U$ such that
$H\in C^2(\overline{\mathbb{B}})$, $J=H^{-1}\in C^2(\overline{U})$, $H(\mathbb{B}^+)=\Omega\cap U$, $H(\mathbb{B}_0)=\p\Omega\cap U$.
Here we denote $\mathbb{B}^+=\mathbb{B}\cap\mathbb{R}^2_+$, $\mathbb{B}^0=\mathbb{B}\cap\p\mathbb{R}^2_+$, and
$\mathbb{R}^2_+=\{(x_1,x_2)\in\mathbb{R}^2: x_2>0\}$. We write $x=H(y)$ and $y=H^{-1}(x)=J(x)$.
Furthermore we can assume (up to a linear transformation) the Jacobian matrix ${\rm Jac}\,H$ satisfies
\be\label{jac}{\rm Jac}\,H(0)=\le(\f{\p H_i}{\p y_j}\ri)_{y=0},\, \f{\p H_i}{\p y_j}(0)=\delta_{ij}=\le\{
\begin{array}{lll}1,\,\,\,&i=j,\\[1.5ex]
0,&i\not=j.\end{array}\ri.\ee
In view of (\ref{E-L-1}),
we have
\be\label{formula}\int_{\Omega\cap U}\nabla u_\epsilon\nabla\varphi dx=\int_{\Omega\cap U}g_\epsilon\varphi dx,\quad
\forall \varphi\in C_0^1({\Omega\cap U}),\ee
where
$$g_\epsilon=\f{1}{\lambda_\epsilon}u_\epsilon e^{(4\pi-\epsilon)u_\epsilon^2}+\alpha u_\epsilon-
  \sum_{i=1}^\ell\sum_{j=1}^{n_i}\f{\gamma_{ij,\epsilon}}{\lambda_\epsilon}e_{ij}.$$
Set
$$\widetilde{u}_\epsilon(y)=u_\epsilon(H(y)),\quad y\in\mathbb{B}^+.$$
Then (\ref{formula}) is transferred to
\be\label{form-2}\sum_{k,\,\ell=1}^2\int_{\mathbb{B}^+}a_{k\ell}\f{\p \widetilde{u}_\epsilon}{\p y_k}\f{\p\psi}{\p y_\ell}dy
=\int_{\mathbb{B}^+}\widetilde{g}_\epsilon\psi dy,\quad \forall \psi\in C_0^1(\mathbb{B}^+),\ee
where
\bna
&&\widetilde{g}_\epsilon=(g_\epsilon\circ H)|{\rm detJac}\, H|,\\[1.2ex]
&&a_{k,\ell}=\sum_{j=1}^2\f{\p J_k}{\p x_j}\f{\p J_\ell}{\p x_j}|{\rm detJac}\, H|,
\ena
and ${\rm detJac}\, H$ denotes the determinant of the Jacobian matrix ${\rm Jac}\,H$.
Note that $a_{k,\ell}\in C^1(\overline{\mathbb{B}^+})$ and that its ellipticity condition is satisfied.

Denote $\widetilde{x}_\epsilon=J(x_\epsilon)=(\widetilde{x}_{1,\epsilon},\widetilde{x}_{2,\epsilon})$
and $\widetilde{x}_\epsilon^\prime=(\widetilde{x}_{1,\epsilon},0)$. Set
$$v_\epsilon(y)=\f{1}{c_\epsilon}\widetilde{u}_\epsilon(\widetilde{x}_\epsilon^\prime+r_\epsilon y),
\quad y\in \overline{\mathbb{B}^+_\epsilon}=\le\{y\in \mathbb{R}^2: \widetilde{x}_\epsilon^\prime+r_\epsilon y\in
\overline{\mathbb{B}}\ri\}.$$
It follows from (\ref{form-2}) that $v_\epsilon$ is a weak solution to the equation
\be\label{form-3}-\f{\p}{\p y_\ell}\le(a_{k\ell}(\widetilde{x}_\epsilon^\prime+r_\epsilon y)\f{\p v_\epsilon}
{\p y_k}(y)\ri)=\f{r_\epsilon^2}{c_\epsilon}\widetilde{g}_\epsilon(\widetilde{x}_\epsilon^\prime+r_\epsilon y),\quad
y\in \mathbb{B}^+_\epsilon.\ee
On one hand, by the definition of $r_\epsilon$ (see (\ref{r-eps-1})), we have
${r_\epsilon^2}{c_\epsilon^{-1}}\widetilde{g}_\epsilon(\widetilde{x}_\epsilon^\prime+r_\epsilon y)$ tends to zero uniformly in
$y\in\mathbb{B}^+_\epsilon$ as $\epsilon\ra 0$. On the other hand we have $|v_\epsilon(y)|\leq 1$ for all $y\in\mathbb{B}^+_\epsilon$.
Note that $\mathbb{B}_\epsilon^+\ra\mathbb{R}^2_+$.
Applying elliptic estimates to (\ref{form-3}) and noticing (\ref{jac}), we obtain $v_\epsilon\ra v$ in $C^1_{\rm loc}
(\overline{\mathbb{R}^2_+})$, where $v$ satisfies
$$\le\{\begin{array}{lll}
-\Delta v=0\quad&{\rm in}\quad &\mathbb{R}^2_+\\[1.2ex]
v=0&{\rm on}&\p\mathbb{R}^2_+\\[1.2ex]
|v|\leq 1& {\rm in} &\mathbb{R}^2_+.
\end{array}\ri.$$
Obviously $v$ can be extended to a bounded weak harmonic function in the whole $\mathbb{R}^2$. Since
$v=0$ on $\p\mathbb{R}^2_+$, Liouville theorem implies that $v\equiv 0$.

We now suppose that there exists some positive number $\nu$ independent of $\epsilon$ such that
\be\label{te-z}r_\epsilon/{\rm dist}(x_\epsilon,\p\Omega)\geq \nu>0.\ee
We can find some constant $C$ depending only on $\nu$ and the bijective map $H$ such that
$$\le|{\widetilde{x}_\epsilon-\widetilde{x}_\epsilon^\prime}\ri|\leq C{r_\epsilon}.$$
Note that
$$v_\epsilon\le(\f{\widetilde{x}_\epsilon-\widetilde{x}_\epsilon^\prime}{r_\epsilon}\ri)=\f{1}{c_\epsilon}
\widetilde{u}_\epsilon(\widetilde{x}_\epsilon)=\f{1}{c_\epsilon}
{u}_\epsilon({x}_\epsilon)=1.$$
We have
$\|v\|_{L^\infty(\mathbb{B}^+_{2C})}=1$, where $\mathbb{B}^+_{2C}=\{y\in \mathbb{R}^2_+:|y|\leq 2C\}$,
since $v_\epsilon\ra v$ in $C_{\rm loc}^1(\overline{\mathbb{R}^2_+})$. This contradicts $v\equiv 0$. Therefore
(\ref{te-z}) is false and our claim (\ref{bb}) follows. \\

  In view of (\ref{bb}), we conclude that
  $$\Omega_\epsilon\ra\mathbb{R}^2\quad{\rm as}\quad \epsilon\ra 0.$$
  Using the argument in Step 2 of the proof of Theorem 1, we have
  \bna&&\psi_\epsilon\ra 1\quad{\rm in}\quad C^1_{\rm loc}(\mathbb{R}^2),\\[1.2ex]
   &&\varphi_\epsilon\ra-\f{1}{4\pi}\log(1+\pi|x|^2)\quad {\rm in}\quad C^1_{\rm loc}(\mathbb{R}^2),\\
   [1.2ex]&&c_\epsilon u_\epsilon\rightharpoonup G\,\,{\rm weakly\,\,in}\,\, W_0^{1,q}(\Omega),\,\,\forall 1<q<2,\\
   [1.2ex]&&c_\epsilon u_\epsilon\ra G\,\,{\rm in}\,\, C^1_{\rm loc}(\overline{\Omega}\setminus\{x_0\}),
   \ena
   where $G$ is a distributional solution to
   $-\Delta G-\alpha G=\delta_{x_0}-\sum_{i=1}^\ell\sum_{j=1}^{n_i}e_{ij}(x_0)e_{ij}$, or equivalently
   \be\label{wso}
   -\int_\Omega G\Delta\varphi dx+\alpha\int_\Omega G\varphi dx=\varphi(x_0)-\sum_{i=1}^\ell\sum_{j=1}^{n_i}e_{ij}(x_0)\int_\Omega
   \varphi e_{ij}dx,\quad \forall \varphi\in C^2(\overline\Omega).
   \ee
   Moreover,
         $$\int_\Omega Ge_{ij}dx=\lim_{\epsilon\ra 0}\int_\Omega c_\epsilon u_\epsilon e_{ij}dx=0,\quad \forall 1\leq j\leq n_i,\,
         1\leq i\leq \ell.$$
         Hence we conclude
         \be\label{G-perp}\int_\Omega Ghdv_g=0,\quad\forall h\in E_\ell.\ee

   If $x_0\in \p\Omega$, testing the equation (\ref{wso}) by $\phi\in C^2(\overline\Omega)\cap W_0^{1,2}(\Omega)$, we have
   $$-\int_\Omega G\Delta\phi dx+\alpha\int_\Omega G\phi dx=0,$$
   since $\phi=0$ on $\p\Omega$ (see \cite{Brezis}, page 288). By the Sobolev embedding theorem, $G\in L^2(\Omega)$.
   By Lemma 6, we have $G\in W^{2,2}(\Omega)\cap W_0^{1,2}(\Omega)$.
   Hence $G$ is an usual weak solution to the equation
   $$
     \left\{\begin{array}{lll}
            -\Delta G-\alpha G=0\quad\rm{in}\quad
            \Omega,\\[1.5ex]
            G\in W_0^{1,2}(\Omega),
         \end{array}\right.
   $$
   and thus $G\equiv 0$ in $\Omega$, since $G\in E_\ell^\perp$ and $0\leq\alpha<\lambda_{\ell+1}(\Omega)$.

   Fourthly, we estimate the supremum (\ref{su-perp})
   under the assumption that $c_\epsilon\ra +\infty$ as
   $\epsilon\ra 0$.
   If $x_0$ lies on the boundary $\p\Omega$,  we set
   $$u_\epsilon^\ast(x)=\le\{\begin{array}{lll}
   &u_\epsilon(x),\quad &x\in \Omega,\\[1.2ex]
   &0,&x\in \mathbb{R}^2\setminus \Omega.\end{array}\ri.$$
   Denote
   $$s_\epsilon=\sup_{\p B_\delta(x_0)}u_\epsilon^\ast,\quad
   \overline{u}_\epsilon^\ast=\le(u_\epsilon^\ast-s_\epsilon\ri)^+,\quad \tau_\epsilon=\int_{B_\delta(x_0)}
   |\nabla\overline{u}_\epsilon^\ast|^2dx.$$
   Since $c_\epsilon u_\epsilon\ra 0$ in $C^1_{\rm loc}(\overline{\Omega}\setminus\{x_0\})\cap L^2(\Omega)$, we have
   \be\label{s-e}s_\epsilon=o_\epsilon(1)c_\epsilon^{-1}\ee
   and
   \bea\nonumber
   \tau_\epsilon&\leq&\int_{B_\delta(x_0)\cap\Omega}|\nabla u_\epsilon|^2dx\\[1.2ex]
   \nonumber&=&1-\int_{\Omega\setminus B_\delta(x_0)}|\nabla u_\epsilon|^2dx+\alpha\int_\Omega u_\epsilon^2dx\\
   [1.2ex]&=&1+\f{o_\epsilon(1)}{c_\epsilon^2}.\label{t-e}
   \eea
   It follows from Lemma 5 that
   \be\label{c-c}\limsup_{\epsilon\ra 0}\int_{B_\delta(x_0)}(e^{4\pi {\overline{u}_\epsilon^\ast}^2/\tau_\epsilon}-1)dx\leq
   \pi\delta^2e.\ee
   In view of (\ref{s-e}) and (\ref{t-e}), there holds on $B_{Rr_\epsilon}(x_\epsilon)$,
   \bna
   (4\pi-\epsilon)u_\epsilon^2&\leq&4\pi {u_\epsilon^\ast}^2\\[1.2ex]
   &\leq&4\pi(\overline{u}_\epsilon^\ast+s_\epsilon)^2\\[1.2ex]
   &=&4\pi{\overline{u}_\epsilon^\ast}^2+8\pi s_\epsilon \overline{u}_\epsilon^\ast+4\pi s_\epsilon^2\\
   [1.2ex]&=&4\pi {\overline{u}_\epsilon^\ast}^2/\tau_\epsilon+o_\epsilon(1).
   \ena
   This together with (\ref{c-c}) leads to
   \bea
   \int_{B_{Rr_\epsilon}(x_\epsilon)}e^{(4\pi-\epsilon)u_\epsilon^2}dx&\leq&\int_{B_{Rr_\epsilon}(x_\epsilon)}
   (e^{4\pi{\overline{u}_\epsilon^\ast}^2}-1)dx+o_\epsilon(1)\nonumber\\[1.2ex]\nonumber&\leq&
   \int_{B_{\delta}(x_0)}
   (e^{4\pi{\overline{u}_\epsilon^\ast}^2}-1)dx+o_\epsilon(1)\\[1.2ex]
   &\leq&\pi \delta^2e+o_\epsilon(1).\label{u-u}
   \eea
   By an analogue of (\ref{u-p-1}), it follows from (\ref{u-u}) that
   $$\sup_{u\in E_\ell^\perp,\,\|u\|_{1,\alpha}\leq 1}\int_\Omega e^{4\pi u^2}dx\leq |\Omega|+\pi\delta^2 e.$$
   Since $\delta>0$ is arbitrary, we get
   $$\sup_{u\in E_\ell^\perp,\,\|u\|_{1,\alpha}\leq 1}\int_\Omega e^{4\pi u^2}dx\leq |\Omega|,$$
   which is impossible. This excludes the possibility of $x_0\in \p\Omega$.

   Now since $x_0\in\Omega$, the Green function $G$ given by (\ref{wso}) can be represented by
          \be\label{Ax0-1}
          G(x)=-\f{1}{2\pi}\log
          |x-x_0|+A_{x_0}+\psi_\alpha(x),\ee where $A_{x_0}$ is a
          constant depending only on $x_0$ and $\alpha$, $\psi_\alpha\in
          C^1(\overline{\Omega})$ and $\psi_\alpha(x_0)=0$. Repeating the argument of deriving (\ref{upperbound}), we get
          \be\label{upb-2}\sup_{u\in E_\ell^\perp,\,\|u\|_{1,\alpha}\leq 1}\int_\Omega e^{4\pi u^2}dx\leq |\Omega|+\pi
          e^{1+4\pi A_{x_0}}.\ee

     Finally we prove the existence of extremal function.
     It suffices to construct a sequence of functions ${\phi}_\epsilon^\ast\in E_\ell^\perp$ with $\|{\phi}_\epsilon^\ast\|_{1,\alpha}=1$
     such that for sufficiently
     small $\epsilon>0$,
     \be\label{gg}
     \int_\Omega e^{4\pi{\phi_\epsilon^\ast}^2}dx>
     |\Omega|+\pi e^{1+4\pi A_{x_0}}.
     \ee
     We shall adapt the test functions constructed in Step 5 of the proof of Theorem 1.
      Let $\phi_\epsilon$ be defined by (\ref{ppp}), $G$ be as in (\ref{wso}), $R=-\log\epsilon$, $c^2$ be as in (\ref{c2-1}),
            and $B$ be as in (\ref{B-1}). In particular $\phi_\epsilon$ satisfies
     the following three properties: $(i)$ $\phi_\epsilon\in W_0^{1,2}(\Omega)$; $(ii)$ $\|\phi_\epsilon\|_{1,\alpha}=1$;
     $(iii)$
      there holds
     $$\int_\Omega e^{4\pi\phi_\epsilon^2}dx\geq
     |\Omega|+\pi e^{1+4\pi A_{x_0}}+4\pi\f{\|G\|_2^2}{c^2}+o(\f{1}{c^2}).
     $$
     Recalling that
     $(e_{ij})$ is a basis of $E_\ell$ verifying (\ref{50}), we set
     $$\widetilde{\phi}_\epsilon=\phi_\epsilon-\sum_{i=1}^\ell\sum_{j=1}^{n_i}(\phi_\epsilon,e_{ij})e_{ij},$$
     where
     $$(\phi_\epsilon,e_{ij})=\int_\Omega \phi_\epsilon e_{ij}dx.$$
     Obviously $\widetilde{\phi}_\epsilon\in E_\ell^\perp$.
     Noting that $e_{ij}\in C^1(\overline{\Omega})$, $R=-\log\epsilon$, $c^2=O(-\log\epsilon)$, $B=O(1)$,
     and $G$ can be represented by (\ref{Ax0-1}), we have
     \bea
     (\phi_\epsilon,e_{ij})&=&\int_{B_{R\epsilon}(x_0)}\le(c+\f{-\f{1}{4\pi}\log(1+\pi\f{r^2}{\epsilon^2})+B}{c}\ri)
      e_{ij}dx\nonumber\\[1.2ex]&&+\int_{B_{2R\epsilon}(x_0)
     \setminus B_{R\epsilon}(x_0)}\f{G-\eta \psi_\alpha}{c} e_{ij}dx+\int_{\Omega\setminus B_{R\epsilon}(x_0)}\f{G}{c} e_{ij}dx\nonumber\\
     [1.2ex]&=& o(\f{1}{\log^2\epsilon}).\label{small}
     \eea
     Here we have used (\ref{G-perp}) to derive
     $$\int_{\Omega\setminus B_{R\epsilon}(x_0)}\f{G}{c} e_{ij}dx=-\int_{B_{R\epsilon}(x_0)}\f{G}{c} e_{ij}dx=O(\epsilon^2(-\log\epsilon)^{5/2})
     =o(\f{1}{\log^2\epsilon}).$$
     By (\ref{small}) and property $(ii)$ of ${\phi}_\epsilon$, we have
     \bea
     &&\widetilde{\phi}_\epsilon={\phi}_\epsilon+o(\f{1}{\log^2\epsilon}),\label{333}\\[1.2ex]
     &&\|\widetilde{\phi}_\epsilon\|_{1,\alpha}^2=1+o(\f{1}{\log^2\epsilon}).\label{444}
     \eea
     Combining (\ref{333}), (\ref{444}) and property $(iii)$ of $\phi_\epsilon$, we obtain
     \bna
     \int_\Omega e^{4\pi\f{\widetilde{\phi}_\epsilon^2}{\|\widetilde{\phi}_\epsilon\|_{1,\alpha}^2}}dx&=&
     \int_\Omega e^{4\pi \phi_\epsilon^2+o(\f{1}{\log\epsilon})}dx\\[1.2ex]
     &\geq&(1+o(\f{1}{\log\epsilon}))\le(|\Omega|+\pi e^{1+4\pi A_{x_0}}+4\pi\f{\|G\|_2^2}{c^2}+o(\f{1}{c^2})\ri)\\
     [1.2ex]&\geq&|\Omega|+\pi e^{1+4\pi A_{x_0}}+4\pi\f{\|G\|_2^2}{c^2}+o(\f{1}{c^2}).
     \ena
     Set $\phi_\epsilon^\ast=\widetilde{\phi}_\epsilon/\|\widetilde{\phi}_\epsilon\|_{1,\alpha}$. Since
     $\widetilde{\phi}_\epsilon\in E_\ell^\perp$, we have $\phi_\epsilon^\ast\in E_\ell^\perp$.
     Moreover $\|\phi_\epsilon^\ast\|_{1,\alpha}=1$ and (\ref{gg}) holds. The contradiction between (\ref{upb-2}) and
     (\ref{gg}) implies that $c_\epsilon$ must be bounded,
     and whence the existence of extremal function follows from (\ref{b-d}) again.
     The proof of Theorem 2 is completely finished. $\hfill\Box$

\section{The Riemannian surface case}

In this section we shall combine Carleson-Chang's result (Lemma 5) and blow-up analysis to
prove Theorems 3 and 4. We follow the lines of \cite{Lijpde,Yang-Tran,Li-Liu-Yang}.
Throughout this section, we denote a geodesic ball centered at $q\in\Sigma$ with radius $r$ by $B_r(q)$,
while a Euclidean ball centered at $x\in\mathbb{R}^2$ with radius $r$ is denoted by $\mathbb{B}_r(x)$.  \\

{\it Proof of Theorem 3.} Let $\alpha$, $0\leq\alpha<\lambda_1(\Sigma)$, be fixed. We divide the proof into several steps.\\

\noindent{\it Step 1. Existence of maximizers for subcritical functionals} \\

In this step, we shall prove for any $0<\epsilon<4\pi$, there exists some $u_\epsilon\in C^1(\Sigma)$ such that
\be\label{bd-1}\|u_\epsilon\|_{1,\alpha}=1,\,\,\int_\Sigma u_\epsilon dv_g=0,\ee
and that
\be\label{sub-c-1}\int_\Sigma e^{(4\pi-\epsilon)u_\epsilon^2}dv_g=\sup_{u\in W^{1,2}(\Sigma),\,
\|u\|_{1,\alpha}\leq1,\,\int_\Sigma u dv_g=0}
\int_\Sigma e^{(4\pi-\epsilon)u^2}dv_g,\ee
where $\|\cdot\|_{1,\alpha}$ is defined as in (\ref{nr-mani}).

 To do this, we choose a maximizing sequence $u_j$ such that $\|u_j\|_{1,\alpha}\leq 1$, $\int_\Sigma u_jdv_g=0$ and
\be\label{sub-2}\int_\Sigma e^{(4\pi-\epsilon)u_j^2}dv_g\ra
\sup_{u\in W^{1,2}(\Sigma),\,
\|u\|_{1,\alpha}\leq1,\,\int_\Sigma u dv_g=0}
\int_\Sigma e^{(4\pi-\epsilon)u^2}dv_g.\ee
It follows from  $0\leq\alpha<\lambda_1(\Sigma)$ that $u_j$ is bounded in $W^{1,2}(\Sigma)$.
Then we can assume, up to a subsequence, $u_j\rightharpoonup u_\epsilon$  weakly  in
  $W^{1,2}(\Sigma)$,  $u_j\ra u_\epsilon$  strongly in  $L^2(\Sigma)$, and
  $u_j\ra u_\epsilon$ a.e.  in $\Sigma$.
  Similarly as in Step 1 of the proof of Theorem 1, we have  $\|u_\epsilon\|_{1,\alpha}\leq 1$  and
  $$\int_\Sigma|\nabla_g u_j-\nabla_g u_\epsilon|^2dv_g\leq
  1-\|u_\epsilon\|_{1,\alpha}^2+o_j(1).$$
  It follows
  from a manifold version of Lions' inequality (\cite{Yang-Tran}, Lemma 3.1) that
  $e^{(4\pi-\epsilon)u_j^2}$ is bounded in $L^q(\Sigma)$ for some $q>1$. Hence $e^{(4\pi-\epsilon)u_j^2}\ra e^{(4\pi-\epsilon)u_\epsilon^2}$
  strongly in $L^1(\Sigma)$. This together with (\ref{sub-2}) leads to (\ref{sub-c-1}).
  Note that $\int_\Sigma u_\epsilon dv_g=0$, since $\int_\Sigma u_j dv_g=0$. We only need to confirm that $\|u_\epsilon\|_{1,\alpha}=1$.
  Suppose not, we have $\|u_\epsilon\|_{1,\alpha}<1$. Set
  ${u^\ast}={u_\epsilon}/\|u_\epsilon\|_{1,\alpha}$.
  Then $u^\ast$ satisfies (\ref{bd-1}) and
  $$\int_\Sigma e^{(4\pi-\epsilon){u^\ast}^2}dv_g>\int_\Sigma e^{(4\pi-\epsilon){u_\epsilon}^2}dv_g,$$
  which contradicts (\ref{sub-c-1}). Therefore $\|u_\epsilon\|_{1,\alpha}=1$.

It is not difficult to check that $u_\epsilon$ satisfies the Euler-Lagrange equation
 \be\label{E-L-2}\le\{
  \begin{array}{lll}
  \Delta_g u_\epsilon-\alpha u_\epsilon=\f{1}{\lambda_\epsilon}u_\epsilon e^{(4\pi-\epsilon)u_\epsilon^2}-\f{\mu_\epsilon}{\lambda_\epsilon}\\[1.5ex]
  \lambda_\epsilon=\int_\Sigma u_\epsilon^2 e^{(4\pi-\epsilon)u_\epsilon^2}dv_g\\
  [1.5ex] \mu_\epsilon=\f{1}{{\rm Vol}_g(\Sigma)}\int_\Sigma u_\epsilon e^{(4\pi-\epsilon)u_\epsilon^2}dv_g,
  \end{array}
  \ri.\ee
  where $\Delta_g$ is the Laplace-Beltrami operator. Applying elliptic estimates to (\ref{E-L-2}), we have that $u_\epsilon\in C^1(\Sigma)$.\\

  \noindent{\it Step 2. Blow-up analysis}\\

  Noting that $$\int_\Sigma e^{(4\pi-\epsilon)u_\epsilon^2}dV_g\leq\int_\Sigma
  \le(1+(4\pi-\epsilon)u_\epsilon^2e^{(4\pi-\epsilon)u_\epsilon^2}\ri)dV_g= {\rm Vol}_g(\Sigma)+(4\pi-\epsilon)\lambda_\epsilon$$
    and $$\lim_{\epsilon\ra 0}\int_\Sigma e^{(4\pi-\epsilon)u_\epsilon^2}dV_g=\sup_{u\in W^{1,2}(\Sigma),\,
\|u\|_{1,\alpha}\leq1,\,\int_\Sigma u dv_g=0}
\int_\Sigma e^{(4\pi-\epsilon)u^2}dv_g,$$
  we have $\liminf_{\epsilon\ra 0}\lambda_\epsilon>0$. It then follows that $\mu_\epsilon/\lambda_\epsilon$ is a bounded sequence.
  Denote $c_\epsilon=|u_\epsilon(x_\epsilon)|=\max_\Sigma |u_\epsilon|$.  If $c_\epsilon$ is bounded, applying elliptic estimates to
  (\ref{E-L-2}), we already conclude the existence of extremal function. Without loss of generality, we assume $x_\epsilon \ra p\in\Sigma$ and
  $c_\epsilon=u_\epsilon(x_\epsilon)\ra +\infty$ as $\epsilon\ra 0$. Take an isothermal coordinate system $(U,\phi)$ near $p$
  such that the metric $g$ can be represented by $g=e^{f}(dx_1^2+dx_2^2)$, where $f\in C^1(\Omega,\mathbb{R})$, $\Omega=\phi(U)\subset\mathbb{R}^2$,
   and $f(0)=0$. Denote $\widetilde{u}_\epsilon=
  u_\epsilon\circ\phi^{-1}$, $\widetilde{x}=\phi^{-1}(x)$ for $x\in\Omega$.

  Let
  $$r_\epsilon={\sqrt{\lambda_\epsilon}}{c_\epsilon^{-1}}e^{-(2\pi-\epsilon/2)c_\epsilon^2},$$
  $$\psi_\epsilon(x)=c_\epsilon^{-1}\widetilde{u}_\epsilon(\widetilde{x}_\epsilon+r_\epsilon
   x),$$ and $$\varphi_\epsilon(x)=c_\epsilon(\widetilde{u}_\epsilon(\widetilde{x}_\epsilon+r_\epsilon
   x)-c_\epsilon)$$
   for $x\in\Omega_\epsilon=\{x\in \mathbb{R}^2:\widetilde{x}_\epsilon+r_\epsilon x\in\Omega\}$.
  By (\ref{E-L-2}), we have
  \be\label{p-s-eq-2}
   -\la_{\mathbb{R}^2} \psi_\epsilon=e^{f(\widetilde{x}_\epsilon+r_\epsilon x)}\le(\alpha r_\epsilon^2\psi_\epsilon+
   c_\epsilon^{-2}\psi_\epsilon e^{(4\pi-\epsilon)(\widetilde{u}_\epsilon^2-c_\epsilon^2)}\ri),
   \ee
   \be\label{phi-eq-2}
   -\la_{\mathbb{R}^2}\varphi_\epsilon=e^{f(\widetilde{x}_\epsilon+r_\epsilon x)}\le(\alpha r_\epsilon^2c_\epsilon^2\psi_\epsilon+\psi_\epsilon e^{(4\pi-\epsilon)(1+\psi_\epsilon)\varphi_\epsilon}\ri),
   \ee
   where $-\Delta_{\mathbb{R}^2}$ denotes the usual Laplacian operator.
   It is easy to see that $\Delta_{\mathbb{R}^2}\psi_\epsilon\ra 0$ in $L^\infty_{\rm loc}(\mathbb{R}^2)$, $|\psi_\epsilon|\leq 1$
   and $\psi_\epsilon(0)=1$. Applying elliptic estimates to (\ref{p-s-eq-2}) and using the Liouville theorem for harmonic function,
   we have
   $$\psi_\epsilon\ra 1\quad {\rm in}\quad C^1_{\rm loc}(\mathbb{R}^2).$$
   Since $\Delta_{\mathbb{R}^2}\varphi_\epsilon$ is bounded in $L^\infty_{\rm loc}(\mathbb{R}^2)$ and
   $\varphi_\epsilon(x)\leq 0=\varphi_\epsilon(0)$ for all $x\in\Omega_\epsilon$, we have by applying elliptic estimates to
   (\ref{phi-eq-2}),
   $$\varphi_\epsilon\ra \varphi=-\f{1}{4\pi}\log(1+\pi|x|^2)\quad {\rm in}\quad C^1_{\rm loc}(\mathbb{R}^2).$$
   Moreover we have
   \be\label{bu-1}\int_{\mathbb{R}^2}e^{8\pi\varphi}dx=1.\ee
   Repeating the argument of proving (\cite{Yang-Tran}, Lemma 4.9), we obtain
   $c_\epsilon u_\epsilon\rightharpoonup G$ weakly in $W^{1,q}(\Sigma)$ for all $1<q<2$, and
   $c_\epsilon u_\epsilon\ra G$ in $C^1_{\rm loc}(\Sigma\setminus\{p\})\cap L^2(\Sigma)$,
   where $G$ is a Green function defined by
  \be\label{Green-1-1}
     \left\{\begin{array}{lll}
            \Delta_g G-\alpha G=\delta_{p}-\frac{1}{{\rm Vol}_g(\Sigma)}\quad\rm{in}\quad
            \Sigma\\[1.5ex]
            \int_\Sigma Gdv_g=0.
         \end{array}\right.\ee
  Clearly $G$ can be represented by
  \be\label{gr-3}G=-\f{1}{2\pi}\log r+A_p+\psi,\ee
  where $r$ denotes the geodesic distance from $p$, $A_p$ is a constant real number, $\psi\in C^1(\Sigma)$ with
  $\psi(p)=0$.\\

  \noindent{\it Step 3. Upper bound estimate}\\

  Similarly as we did in Step 4 of the proof of Theorem 1, we obtain by using Carleson-Chang's result (Lemma 5)
  \be\label{sha}\lim_{R\ra+\infty}\limsup_{\epsilon\ra 0}\int_{B_{Rr_\epsilon}(x_\epsilon)} e^{(4\pi-\epsilon) u_\epsilon^2}dv_g\leq \pi e^{1+4\pi A_p},\ee
   where $A_p$ is given by (\ref{gr-3}). Note that
   \bna
   \int_{B_{Rr_\epsilon}(x_\epsilon)}e^{(4\pi-\epsilon)u_\epsilon}dv_g&=&(1+o_\epsilon(1))
   \int_{\mathbb{B}_{Rr_\epsilon}(\widetilde{x}_\epsilon)}e^{(4\pi-\epsilon)\widetilde{u}_\epsilon}dx\\
   [1.2ex]&=&(1+o_\epsilon(1))
   \int_{\mathbb{B}_{R}(0)}e^{(4\pi-\epsilon)\widetilde{u}_\epsilon}r_\epsilon^2dx\\
   &=&(1+o_\epsilon(1))\f{\lambda_\epsilon}{c_\epsilon^2}\int_{\mathbb{B}_{R}(0)}e^{8\pi\varphi}dx.
   \ena
   This together with (\ref{bu-1}) implies
   $$\lim_{R\ra+\infty}\limsup_{\epsilon\ra 0}\int_{B_{Rr_\epsilon}(x_\epsilon)} e^{(4\pi-\epsilon) u_\epsilon^2}dv_g=
   \limsup_{\epsilon\ra 0}\f{\lambda_\epsilon}{c_\epsilon^2},$$
   which together with (\ref{sha}) and an analogue of (\cite{Yang-Tran}, Lemma 4.6) leads to
   \bea \label{ub-1}\sup_{u\in W^{1,2}(\Sigma),\,\|u\|_{1,\alpha}\leq 1,\,\int_\Sigma u dv_g=0}
\int_\Sigma e^{4\pi u^2}dv_g=\limsup_{\epsilon\ra 0}\int_\Sigma e^{(4\pi-\epsilon) u_\epsilon^2}dv_g\leq {\rm Vol}_g(\Sigma)+\pi e^{1+4\pi A_p}.\eea

     \noindent  {\it Step 4. Existence of extremal function}\\

   In this step we will construct a blow-up sequence $\phi_\epsilon$ such that
   \be\label{nrm}\int_\Sigma|\nabla_g\phi_\epsilon|^2dv_g-\alpha\int_\Sigma (\phi-\overline{\phi}_\epsilon)^2dv_g=1\ee
   and
   \be\label{gr}\int_\Sigma e^{4\pi (\phi_\epsilon-\overline{\phi}_\epsilon)^2}dv_g>{\rm Vol}(\Sigma)+\pi e^{1+4\pi A_p}\ee
   for sufficiently small $\epsilon>0$, where
   $$\overline{\phi}_\epsilon=\f{1}{{\rm Vol}_g(\Sigma)}\int_\Sigma \phi_\epsilon dv_g.$$
   The contradiction between (\ref{gr}) and (\ref{ub-1}) implies that $c_\epsilon$ must be bounded and elliptic estimates
   imply the existence of the desired extremal function. This completes the proof of Theorem 3.

   Now we construct $\phi_\epsilon$ verifying (\ref{nrm}) and (\ref{gr}).
   Note that the Green function $G$ defined as in (\ref{Green-1-1}) has the representation (\ref{gr-3}).
   Set
     \be\label{test}\phi_\epsilon=\le\{
     \begin{array}{llll}
     &c+\f{-\f{1}{4\pi}\log(1+\pi\f{r^2}{\epsilon^2})+B}{c}
     \quad &{\rm for} &r\leq R\epsilon\\[1.5ex]
     &\f{G-\eta \psi}{c}\quad &{\rm for} & R\epsilon<
     r<2R\epsilon\\[1.2ex]
     &\f{G}{c}\quad &{\rm for} & r\geq 2R\epsilon,
     \end{array}
     \ri.\ee
     where $R=-\log\epsilon$, $\eta\in C_0^\infty(B_{2R\epsilon}(p))$ verifying that $\eta=1$ on $B_{R\epsilon}(p)$ and
     $\|\nabla_g \eta\|_{L^\infty}
     =O(\f{1}{R\epsilon})$, $B$ is a constant to be determined
     later, and $c$ depending only on $\epsilon$ will also be
     chosen later such that $R\epsilon\ra 0$ and $R\ra +\infty$.
     In order to assure that $\phi_\epsilon\in W^{1,2}(\Sigma),$ we set
     $$
     c+\f{1}{c}\le(-\f{1}{4\pi}\log(1+\pi R^2)+B\ri)
     =\f{1}{c}\le(-\f{1}{2\pi}\log (R\epsilon)+A_{p}\ri),
     $$
     which gives
     \be\label{2pic2}
     2\pi c^2=-\log\epsilon-2\pi B+2\pi A_{p}+\f{1}{2}\log \pi
     +O(\f{1}{R^2}).
     \ee
     We calculate
     \bna
     \int_{\Sigma}|\nabla_g \phi_\epsilon|^2dv_g&=&\f{1}{4\pi c^2}\le(
     2\log\f{1}{\epsilon}+\log\pi-1+4\pi A_{p}+4\pi\alpha\|G\|_2^2\ri.\\
     &&\le.+O(\f{1}{R^2})+
     O(R\epsilon\log(R\epsilon))\ri),
     \ena
     \bna
     \int_\Sigma \phi_\epsilon dv_g&=&\f{1}{c}\le(\int_{r\geq 2R\epsilon}Gdv_g+O(R\epsilon\log(R\epsilon))\ri)\\
     &=&\f{1}{c}\le(-\int_{r< 2R\epsilon}Gdv_g+O(R\epsilon\log(R\epsilon))\ri)\\
     &=&\f{1}{c}O(R\epsilon\log(R\epsilon)),
     \ena
     and
     \bna
     \int_\Sigma(\phi_\epsilon-\overline{\phi}_\epsilon)^2dv_g&=&\int_\Sigma\phi_\epsilon^2dv_g-\overline{\phi}_\epsilon^2
     {\rm Vol}_g(\Sigma)\\&=&\f{1}{c^2}\le(\int_\Sigma G^2dv_g+O(R\epsilon\log(R\epsilon))\ri).
     \ena
     This yields
     \bna \|\phi_\epsilon-\overline\phi_\epsilon\|_{1,\alpha}^2&=&\int_\Sigma|\nabla_g\phi_\epsilon|^2dv_g-\alpha
     \int_\Sigma(\phi_\epsilon-\overline{\phi}_\epsilon)^2dv_g\\
     &=&
     \f{1}{4\pi c^2}\le(
     2\log\f{1}{\epsilon}+\log\pi-1+4\pi A_{p}
     +O(\f{1}{R^2})+
     O(R\epsilon\log(R\epsilon))\ri)\ena
     Let $\phi_\epsilon$ satisfy (\ref{nrm}), i.e. $\|\phi_\epsilon-\overline\phi_\epsilon\|_{1,\alpha}=1$. Then we have
     \be\label{c2}
     c^2=-\f{\log\epsilon}{2\pi}+\f{\log\pi}{4\pi}-\f{1}{4\pi}+A_{p}
     +O(\f{1}{R^2})+O(R\epsilon\log(R\epsilon)).
     \ee
     It follows from (\ref{2pic2}) and (\ref{c2}) that
     \be{\label{B}}
     B=\f{1}{4\pi}+O(\f{1}{R^2})+O(R\epsilon\log(R\epsilon)).
     \ee
     Clearly we have on $B_{R\epsilon}(p)$
     $$4\pi(\phi_\epsilon-\overline{\phi}_\epsilon)^2\geq 4\pi
     c^2-2\log(1+\pi\f{r^2}{\epsilon^2})+8\pi B.$$
     This together with (\ref{c2}) and (\ref{B}) yields
     \be\label{BRE}
     \int_{B_{R\epsilon}(p)} e^{4\pi(\phi_\epsilon-\overline{\phi}_\epsilon)^2}dv_g\geq
     \pi e^{1+4\pi A_{p}}
     +O(\f{1}{(\log\epsilon)^2}).
     \ee
     On the other hand,
     \bea
      \label{O-BRE}\int_{\Sigma\setminus
      B_{R\epsilon}(p)}e^{4\pi(\phi_\epsilon-\overline{\phi}_\epsilon)^2}dv_g&\geq&\int_{\Sigma\setminus
      B_{2R\epsilon}(p)}(1+4\pi\phi_\epsilon^2)dv_g{\nonumber}\\
      &\geq& {\rm Vol}_g(\Sigma)+4\pi\f{\|G\|_2^2}{c^2}+o(\f{1}{c^2}).
     \eea
     Recalling (\ref{c2}) and combining (\ref{BRE}) and (\ref{O-BRE}),
     we conclude (\ref{gr}) for sufficiently small $\epsilon>0$.
     $\hfill\Box$\\

     {\it Proof of Theorem 4.}
     Let $\lambda_1(\Sigma)<\lambda_2(\Sigma)<\cdots$ be all distinct eigenvalues of the Laplace-Beltrami operator
     $\Delta_g$, and $E_{\lambda_i(\Sigma)}$'s be
     associated eigenfunction spaces. It is known that $\lambda_i(\Sigma)\ra +\infty$ as $i\ra+\infty$
     and each space $E_{\lambda_i(\Sigma)}$ has finite dimension (see \cite{Chavel}, Chapter I, Page 8).
     We can assume
     $${\rm dim}E_{\lambda_i(\Sigma)}=n_i,\quad i=1,2,\cdots.$$
     Take a basis $(e_{ij})$ $(1\leq j\leq n_i,1\leq i\leq\ell)$ of $E_\ell$ verifying
     \bna
     &&E_{\lambda_i(\Sigma)}={\rm span}\{e_{i1},\cdots,e_{in_i}\},\quad i=1,\cdots,\ell,\\
     [1.2ex]&&E_\ell={\rm span}\{e_{11},\cdots,e_{1n_1}, e_{21},\cdots,e_{2,n_2},\cdots,e_{\ell 1},\cdots,e_{\ell n_\ell}\},\\
     [1.2ex] &&\int_\Sigma|e_{ij}|^2dv_g=1,\\
     [1.2ex] &&\int_\Sigma e_{ij}e_{kl}dv_g=0,\,\, i\not= k\,\,{\rm or}\,\, j\not=l.
     \ena
     Similar to Step 1 of the proof of Theorem 3, for any $\epsilon$, $0<\epsilon<4\pi$, there exists
     some $u_\epsilon\in E_\ell^\perp\cap C^1(\Sigma)$ with
 $\|u_\epsilon\|_{1,\alpha}=1$ such that
 $$\label{subcrit1}\int_\Sigma e^{(4\pi-\epsilon)u_\epsilon^2}dv_g=
 \sup_{u\in E_\ell^\perp,\,\|u\|_{1,\alpha}\leq 1}\int_\Sigma e^{(4\pi-\epsilon)u^2}dv_g.$$
 Moreover $u_\epsilon$ satisfies the
 Euler-Lagrange equation
  $$\label{Eule}\le\{
  \begin{array}{lll}
  \Delta_g u_\epsilon-\alpha u_\epsilon=\f{1}{\lambda_\epsilon}u_\epsilon e^{(4\pi-\epsilon)u_\epsilon^2}-
  \f{\mu_\epsilon}{\lambda_\epsilon}-\sum_{i=1}^\ell\sum_{j=1}^{n_i}\f{\gamma_{ij,\epsilon}}{\lambda_\epsilon}e_{ij}\,\,\,{\rm in}\,\,\,
  \Sigma,\\[1.2ex] u_\epsilon\in E_\ell^\perp\cap C^1(\Sigma),\\[1.2ex]
  \lambda_\epsilon=\int_\Sigma u_\epsilon^2 e^{(4\pi-\epsilon)u_\epsilon^2}dv_g,\\[1.2ex]
  \mu_\epsilon=\f{1}{{\rm Vol}_g(\Sigma)}\int_\Sigma u_\epsilon e^{(4\pi-\epsilon)u_\epsilon^2}dv_g,\\[1.2ex]
  \gamma_{ij,\epsilon}=\int_\Sigma e_{ij}u_\epsilon e^{(4\pi-\epsilon)u_\epsilon^2}dv_g.
  \end{array}
  \ri.$$
  Let $c_\epsilon=\max_\Sigma |u_\epsilon|$. Without loss of generality, we assume
  $c_\epsilon=u_\epsilon(x_\epsilon)\ra +\infty$  and $x_\epsilon\ra p\in\Sigma$ as $\epsilon\ra 0$.
   Take an isothermal coordinate system $(U,\phi)$ near $p$. Denote $\widetilde{u}_\epsilon=
  u_\epsilon\circ\phi^{-1}$, $\widetilde{x}=\phi^{-1}(x)$ for $x\in\phi(U)\subset\mathbb{R}^2$.
  Perform the same blow-up analysis as in the proof of Theorem 3. There holds
   \bna&&\varphi_\epsilon\ra \varphi=-\f{1}{4\pi}\log(1+\pi|x|^2)\quad {\rm in}\quad C^1_{\rm loc}(\mathbb{R}^2),\\
   [1.2ex]&&c_\epsilon u_\epsilon\rightharpoonup G\,\,{\rm weakly\,\,in}\,\, W^{1,q}(\Sigma),\,\,\forall 1<q<2,\\
   [1.2ex]&&c_\epsilon u_\epsilon\ra G\,\,{\rm in}\,\, C^1_{\rm loc}(\Sigma\setminus\{p\})\cap L^2(\Sigma),\ena
   where $\varphi_\epsilon(x)=c_\epsilon(\widetilde{u}_\epsilon(\widetilde{x}_\epsilon+r_\epsilon x)-c_\epsilon)$,
   $G$ is a Green function defined by
   \be\label{gn}\Delta_g G-\alpha G=\delta_{p}-\f{1}{{\rm Vol}_g(\Sigma)}-\sum_{i=1}^\ell\sum_{j=1}^{n_i}e_{ij}(p)e_{ij}.\ee
   Since $u_\epsilon\in E_\ell^\perp$, we have
   $$\int_\Sigma Ge_{ij}dv_g=\lim_{\epsilon\ra 0}\int_\Sigma c_\epsilon u_\epsilon e_{ij}dv_g=0,\quad
   \forall 1\leq j\leq n_i,\,1\leq i\leq\ell.$$
  Clearly $G$ can be written as
  \be\label{gr-32}G=-\f{1}{2\pi}\log r+A_p+\psi,\ee
  where $r$ denotes the geodesic distance from $p$, $A_p$ is a constant real number, $\psi\in C^1(\Sigma)$ with
  $\psi(p)=0$. Using Carleson-Chang's result (Lemma 5), we obtain
  \bea \label{ub-12}\sup_{u\in E_\ell^\perp,\,\|u\|_{1,\alpha}\leq 1,\,\int_\Sigma u dv_g=0}
\int_\Sigma e^{4\pi u^2}dv_g=\limsup_{\epsilon\ra 0}\int_\Sigma e^{(4\pi-\epsilon) u_\epsilon^2}dv_g\leq {\rm Vol}_g(\Sigma)+\pi e^{1+4\pi A_p}.\eea

   Now we will construct a sequence of functions $\phi_\epsilon^\ast$ such that $\phi_\epsilon^\ast\in E_\ell^\perp$, $\int_\Sigma
   \phi_\epsilon^\ast dv_g=0$ and
   \be\label{gr-3-1}\int_\Sigma e^{4\pi {\phi_\epsilon^\ast}^2}dv_g>{\rm Vol}_g(\Sigma)+\pi e^{1+4\pi A_p}\ee
   for sufficiently small $\epsilon>0$.
   The contradiction between (\ref{gr-3-1}) and (\ref{ub-12}) implies that $c_\epsilon$ must be bounded and elliptic estimates
   lead to the existence of the desired extremal function. This completes the proof of Theorem 4.

     Let $\phi_\epsilon$ be defined by (\ref{test}), $G$ be as in (\ref{gn}), $R=-\log\epsilon$, $c^2$ be as in (\ref{c2}),
            and $B$ be as in (\ref{B}). In particular $\phi_\epsilon$ satisfies
      \be\label{n-1-2}\int_\Sigma|\nabla_g\phi_\epsilon|^2dv_g-\alpha\int_\Sigma \phi_\epsilon^2dv_g=1\ee
      and
     \be\label{ggg-1}\int_\Sigma e^{4\pi(\phi_\epsilon-\overline{\phi}_\epsilon)^2}dv_g\geq
     {\rm Vol}_g(\Sigma)+\pi e^{1+4\pi A_{p}}+4\pi\f{\|G\|_2^2}{c^2}+o(\f{1}{c^2}),
     \ee
     where
     $$\overline{\phi}_\epsilon=\f{1}{{\rm                                                                                                                   Vol}_g(\Sigma)}\int_\Sigma \phi_\epsilon dv_g.$$
      Set
     $$\widetilde{\phi}_\epsilon=\phi_\epsilon-\overline{\phi}_\epsilon-\sum_{i=1}^\ell\sum_{j=1}^{n_i}(\phi_\epsilon
     -\overline{\phi}_\epsilon,e_{ij})e_{ij},$$
     where
     $$(\phi_\epsilon-\overline{\phi}_\epsilon,e_{ij})=\int_\Sigma (\phi_\epsilon-\overline{\phi}_\epsilon) e_{ij}dv_g.$$
     Obviously $\widetilde{\phi}_\epsilon\in E_\ell^\perp$.
     Note that $e_{ij}\in C^1({\Sigma})$, $R=-\log\epsilon$, $c^2=O(-\log\epsilon)$, $B=O(1)$,
     and $G$ can be represented by (\ref{gr-32}). We calculate
     \bna
     (\phi_\epsilon-\overline{\phi}_\epsilon,e_{ij})&=&\int_{B_{R\epsilon}(p)}\le(c+\f{-\f{1}{4\pi}\log(1+\pi\f{r^2}{\epsilon^2})+B}{c}-
     \overline{\phi}_\epsilon\ri)
      e_{ij}dv_g\nonumber\\[1.2ex]&&+\int_{B_{2R\epsilon}(p)
     \setminus B_{R\epsilon}(p)}\le(\f{G-\eta \psi_\alpha}{c}-\overline{\phi}_\epsilon\ri) e_{ij}dv_g+
     \int_{\Sigma\setminus B_{R\epsilon}(p)}\le(\f{G}{c}-\overline{\phi}_\epsilon\ri) e_{ij}dv_g\nonumber\\
     [1.2ex]&=& I+II+III.
     \ena
     Since $G\in E_\ell^\perp$, we have
     $$\int_{\Sigma\setminus B_{R\epsilon}(p)}\f{G}{c}dv_g=-\int_{B_{R\epsilon}(p)}\f{G}{c}dv_g=\f{1}{c}O(R^2\epsilon^2\log(R\epsilon)).$$
     Note also that $\overline{\phi}_\epsilon=\f{1}{c}O(R\epsilon\log(R\epsilon))$.
     Hence
     $$\label{III}III=\f{1}{c}O(R\epsilon\log(R\epsilon)).$$
     Clearly
     $$\label{II}I=O(cR^2\epsilon^2),\quad II=\f{1}{c}O(R^2\epsilon^2\log(R\epsilon)).$$
     Therefore
     $$\label{kk}(\phi_\epsilon-\overline{\phi}_\epsilon,e_{ij})=O(R^2\epsilon^2\sqrt{-\log{\epsilon}})=
     o(\f{1}{\log^2\epsilon}).$$
     This together with (\ref{n-1-2}) leads to
     \bea
     &&\widetilde{\phi}_\epsilon={\phi}_\epsilon-\overline{\phi}_\epsilon+o(\f{1}{\log^2\epsilon}),\label{one}\\[1.2ex]
     &&\|\widetilde{\phi}_\epsilon\|_{1,\alpha}^2=1+o(\f{1}{\log^2\epsilon}).\label{two}
     \eea
     Combining (\ref{one}), (\ref{two}) and (\ref{ggg-1}), we obtain
     \bna
     \int_\Sigma e^{4\pi\f{\widetilde{\phi}_\epsilon^2}{\|\widetilde{\phi}_\epsilon\|_{1,\alpha}^2}}dv_g&=&
     \int_\Sigma e^{4\pi (\phi_\epsilon-\overline{\phi}_\epsilon)^2+o(\f{1}{\log\epsilon})}dv_g\\[1.2ex]
     &\geq&(1+o(\f{1}{\log\epsilon}))\le({\rm Vol}_g(\Sigma)+\pi e^{1+4\pi A_{p}}+4\pi\f{\|G\|_2^2}{c^2}+o(\f{1}{c^2})\ri)\\
     [1.2ex]&\geq&{\rm Vol}_g(\Sigma)+\pi e^{1+4\pi A_{p}}+4\pi\f{\|G\|_2^2}{c^2}+o(\f{1}{c^2}).
     \ena
     Set $\phi_\epsilon^\ast=\widetilde{\phi}_\epsilon/\|\widetilde{\phi}_\epsilon\|_{1,\alpha}$. Since
     $\widetilde{\phi}_\epsilon\in E_\ell^\perp$, we have $\phi_\epsilon^\ast\in E_\ell^\perp$.
     Moreover $\|\phi_\epsilon^\ast\|_{1,\alpha}=1$. Since $\Delta_g e_{ij}=\lambda_i e_{ij}$, we have
     $\int_\Sigma e_{ij}dv_g=0$, and whence $\int_\Sigma{\phi}_\epsilon^\ast dv_g=0$.
     Therefore $\phi_\epsilon^\ast$ is the desired function sequence verifying (\ref{gr-3-1}). $\hfill\Box$\\


 {\bf Acknowledgements.} This work is supported by the National Science Foundation of China (Grant No.11171347 and Grant
 No. 11471014).
 The author also gratefully acknowledges the helpful comments and suggestions of the referees, which have improved
 the presentation.

\end{document}